\documentclass[preprint,10pt,nopreprintline]{elsarticle}

\usepackage[english]{babel}
\usepackage{geometry}
\usepackage{amsmath,amsfonts,amssymb,amsthm}
\usepackage[final]{hyperref}
\usepackage{graphicx}
\usepackage{mathtools}
\usepackage{physics}
\usepackage{mathtools}
\usepackage[dvipsnames]{xcolor}
\usepackage{pdflscape}
\usepackage{etoolbox}
\usepackage{fancyhdr}
\usepackage{lastpage}
\usepackage{ifdraft}
\usepackage{hyphenat}
\usepackage{natbib}
\ifdraft{\usepackage{showlabels}}{}
\usepackage{enumitem}
\usepackage{siunitx}
\usepackage{booktabs}
\usepackage{multirow}
\usepackage{array}
\usepackage{caption}
\usepackage{float}
\usepackage{rotating}

\DeclareMathOperator*{\argmin}{\arg\!\min}

\usepackage{lineno}

\begin{document}

\begin{frontmatter}



\title{Sparse Identification of Nonlinear Dynamics for Stochastic Delay Differential Equations} 


\author[1]{Dimitri Breda}
\ead{dimitri.breda@uniud.it}
\author[2]{Dajana Conte}
\ead{dajconte@unisa.it}
\author[3]{Raffaele D'Ambrosio}
\ead{raffaele.dambrosio@univaq.it}
\author[2]{Ida Santaniello}
\ead{isantaniello@unisa.it}
\author[1]{Muhammad Tanveer}
\ead{tanveer.muhammad@spes.uniud.it}

\affiliation[1]{organization={CDLab – Computational Dynamics Laboratory Department of Mathematics, Computer Science and Physics – University of Udine},
                addressline={via delle Scienze 206}, 
                city={Udine},
                postcode={33100}, 
                country={Italy}}
\affiliation[2]{organization={Department of Mathematics, University of Salerno},
                addressline={Via Giovanni Paolo II 132}, 
                city={Fisciano (SA)},
                postcode={84084}, 
                country={Italy}}
\affiliation[3]{organization={Department of Information Engineering and Computer Science and Mathematics, University of L’Aquila},
                addressline={Vetoio, Loc. Coppito}, 
                city={L’Aquila},
                postcode={67100}, 
                country={Italy}}

\begin{abstract}
A general framework for recovering drift and diffusion dynamics from sampled trajectories is presented for the first time for stochastic delay differential equations. The core relies on the well-established SINDy algorithm for the sparse identification of nonlinear dynamics. The proposed methodology combines recently proposed high-order estimates of drift and covariance for dealing with stochastic problems with augmented libraries to handle delayed arguments. Three different strategies are discussed in view of exploiting only realistically available data. A thorough comparative numerical investigation is performed on different models, which helps guiding the choice of effective and possibly outperforming schemes.
\end{abstract}


%
%
%
%

\end{frontmatter}



\section{Introduction} \label{s_{i}ntroduction}
{Stochastic differential equations (SDEs) provide a natural framework for modeling systems subject to random fluctuations. Among them, stochastic ordinary differential equations (SODEs) describe the evolution of systems whose future state depends on their current state and random perturbations, typically modeled by Brownian motion. In contrast, stochastic delay differential equations (SDDEs) incorporate memory effects by allowing the system's dynamics to depend not only on the present but also on past states. The analysis of such systems have been widely explored in the literature, both from the theoretical and numerical perspectives, as well as from a data-driven point of view. \cite{buckwar2000introduction,higham2001algorithmic, kloeden1992stochastic, mohammed1998stochastic,oksendal2003stochastic} provide an excellent review of these and relevant topics. 

In the last years, there has been a growing interest in the identification of stochastic dynamical systems from observed time-series data. In this context, Sparse Identification of Nonlinear Dynamics (SINDy) has emerged as a powerful data-driven algorithm for discovering governing equations from trajectories \cite{bk19, bpk16}. The principal aim of SINDy is to identify parsimonious and interpretable models from data, although its performance may decline when only limited measurements are available, due to its sensitivity to data sparsity and noise. Several works have applied SINDy to SODEs \cite{boninsegna2018sparse, wanner2024higher} and deterministic delay differential equations (DDEs) \cite{bbt24,kopeczi23,pec24,sandoz2023sindy}, while its extension to SDDEs has not yet been fully investigated, owing to its technical and computational difficulties. This work aims to address these challenges by extending the SINDy framework to the SDDEs setting. In this regard, we contrast our approach with \cite{han2024approximation}, which is restricted to the Stratonovich formulation and whose reliance on small-delay Taylor approximations limits generalization to arbitrary delays.

We explore and compare several estimations strategies, including Kramers–Moyal (KM) formulas, forward differences (FD), central differences (CD) and trapezoidal-like schemes (TR), all based on Itô–Taylor expansions, for reconstructing drift and diffusion from time-series data and extend these estimation strategies to the framework of SDDEs, relying on an augmented state with delay. The resulting approximations are then used as input for the sparse regression performed by the SINDy algorithm. In this respect, we investigate three practical approaches (A, B1, B2) for recovering drift and diffusion data using the SINDy algorithm combined with the discussed Itô–Taylor-based estimators, under different data availability scenarios. The proposed approaches address the fundamental challenge of obtaining reliable estimates at each time step and are designed to complement the developed estimation methods. These strategies differ in terms of data requirements, statistical robustness, obtained accuracy and computational cost in terms of CPU time.

The paper is organized as follows.
In Section \ref{sec:theor_back}, we introduce the notation and main concepts related to SODEs and SDDEs, focusing on how to recover the KM formulas and their higher-order extensions for estimating drift and diffusion terms, based on Itô–Taylor expansions.
Section \ref{sec:sindy} provides a concise overview of the original SINDy framework for ordinary differential equations (ODEs) and its extensions first to delay differential equations (DDEs) and second to SODEs. The discussion then moves to SDDEs, highlighting how the method can be adapted to account for both delayed and stochastic components.
In Section \ref{sec:estimates}, we turn to the practical use of the estimators introduced in Section \ref{sec:theor_back} depending on the available data. We propose three different strategies for recovering drift and diffusion terms: Approach A, which generates synthetic paths from a single observed trajectory; Approach B1, which computes averaged estimates before sparse regression using multiple trajectories; and Approach B2, which applies sparse regression separately to each trajectory and averages the resulting models. 

Their comparative performance is evaluated in Section \ref{sec:tests}, which presents a computational analysis on three representative SDDE models. For each model, the three approaches A, B1 and B2 are tested in combination with all four estimation methods KM, FD, CD and TR  discussed in Section \ref{sec:theor_back}. Synthetic data are generated via Euler–Maruyama simulations. Reconstruction quality is assessed using absolute errors on the identified sparse coefficients and root mean square errors (RMSE) on both reconstructed trajectories and estimated drift and diffusion components. Additionally, we analyze the impact of data availability and CPU time on the overall performance.

\section{Theoretical background}\label{sec:theor_back}
In Section \ref{sec:sde} we recall the basics of SODEs first and SDDEs then: notation, key ingredients and well-posedness of the relevant initial value problems. In Section \ref{sec:ito_taylor} we briefly review the associated It\^{o}-Taylor expansions, on which Kramers-Moyal and higher-order estimates of drift and diffusion erly. They form the necessary theoretical foundation for extending the sparse identification framework of SINDy to both SODEs and SDDEs. Useful general references for the underlying theoretical framework are \citep{hk21,higham2001algorithmic,kloeden1992stochastic,mao1996robustness,mohammed1998stochastic,oksendal2003stochastic}.
\subsection{Stochastic ordinary and delay differential equations}\label{sec:sde}
SODEs describe the evolution in time $t\geq0$ of a real random variable $X(t)\in\mathbb{R}^n$ through the interplay of deterministic and stochastic forces. In its classical formulation without delays a SODE reads
\begin{equation}\label{eq:sode}
\dd X(t)=f(X(t))\dd t+g(X(t))\dd W(t),
\end{equation}
where $W(t)$ is a $q$-dimensional standard Wiener process, the \emph{drift} function $f:\mathbb{R}^n\to \mathbb{R}^n$ governs the deterministic evolution and represents the expected direction of motion at each point, while the \emph{diffusion} function $g: \mathbb{R}^n\to \mathbb{R}^{n\times q}$ controls the intensity and structure of random fluctuations. We assume the \emph{covariance} matrix $G\coloneqq gg^{T}/2\in\mathbb{R}^{n\times n}$ to be positive definite, thus ensuring that $g$ is non-degenerate and can be uniquely determined from $G$ up to a right orthogonal transformation \citep{o2021interacting,raj2023efficient}. The latter assumption is not restrictive and indeed many physical systems affected by noise fit into this framework  \citep{gardiner1985handbook,oksendal2003stochastic}. Of course \eqref{eq:sode} is just a symbolic notation for
\begin{equation*}
X(t)=X(0)+\int_{0}^{t}f(X(\sigma))\dd\sigma+\int_{0}^{t}g(X(\sigma))\dd W(\sigma),
\end{equation*}
where the last integral is intended in the It\^{o} sense as far as this work is concerned. Under the standard global Lipschitz and linear growth conditions
\begin{align*}
\|f(x)-f(y)\| + \|g(x)-g(y)\| &\leq L\|x-y\|, \quad x,y \in \mathbb{R}^{n},\\
\|f(x)\| + \|g(x)\| &\leq C(1+\|x\|), \quad x \in \mathbb{R}^{n},
\end{align*}
for $L$ and $C$ positive constants, the initial value problem for \eqref{eq:sode} admits a unique strong solution with continuous sample paths (for a proof see, e.g., \cite{oksendal2003stochastic}). We observe that alternative hypotheses may also be assumed, e.g., one-sided Lipschitz continuous drift and locally Lipschitz continuous diffusion, in order to ensure bounded moments for the solution \cite{highamWeak}.

\bigskip
Many real-world systems exhibit memory effects, for which the current evolution of the state depends also on past values \citep{da2004color,van1992stochastic}. In the presence of randomness SDDEs capture this phenomenon by incorporating explicit time delays into the dynamics \citep{buckwar2000introduction,guillouzic1999small,mao1996robustness,tian2007stochastic}. Here we consider SDDEs with a single discrete and constant delay $\tau>0$ affecting both drift and diffusion:
\begin{equation}\label{eq:sdde}
\dd X(t) = f(X(t),X({t-\tau}))\dd t + g(X(t),X({t-\tau}))\dd W({t}).
\end{equation}
Extension to multiple delays is not conceptually demanding, even though it might involve some technicalities. In \eqref{eq:sdde} drift and diffusion are described, respectively, by $f: \mathbb{R}^n\times \mathbb{R}^n\to \mathbb{R}^n$ and $g: \mathbb{R}^n\times \mathbb{R}^n\to \mathbb{R}^{n\times q}$. A related initial value problem requires specifying an initial history through a random variable $\varphi(t)$ defined on $[-\tau,0]$ by imposing $X({t})=\varphi(t)$ for $t\in[-\tau,0]$. Unlike standard SODEs, solutions to SDDEs are generally non-Markovian due to their dependence on the past \cite{mao1996robustness}. Nevertheless, the initial value problem for \eqref{eq:sdde} becomes well-posed under appropriate Lipschitz and linear growth conditions whose general form reads
\begin{align*}
\|f(\varphi) - f(\psi)\| + \|g(\varphi) - g(\psi)\| \leq L \|\varphi - \psi\|_C, \quad \varphi,\psi \in C,\\
\|f(\varphi)\| + \|g(\varphi)\| \leq C(1 + \|\varphi\|_C), \quad \varphi \in C,
\end{align*}
for drift and diffusion given as maps $C\to\mathbb{R}^{n}$ acting on the full history defined as a function $[-\tau,0]\ni s\mapsto X(t)(s)\coloneqq X(t+s)$ for $X(t)\in C\coloneqq C([-\tau, 0]; \mathbb{R}^n)$ the classic state space of continuous functions equipped with the uniform norm $\|\varphi\|_C\coloneqq \max_{s\in[-\tau,0]} \|\varphi(s)\|$ \citep{buckwar2000introduction,mao1996robustness}.

As anticipated qualitatively, our goal is to estimate the assumed to be unknown drift $f$ and diffusion $g$ from trajectory data given in terms of the $m$ time samples $X(t_{1}),\ldots,X(t_{m})$, possibly related to multiple paths. Note once more that one can recover $g$ uniquely once $G$ has been identified under the assumption that the latter is positive definite.
\subsection{Estimates of drift and diffusion via It\^{o}-Taylor Expansions}\label{sec:ito_taylor}
Practical estimation methods for $f$ and $G$ are based on It\^{o}-Taylor expansions, which express increments of test functions of a stochastic process in terms of differential operators and multiple stochastic integrals. They provide a systematic way to approximate solutions of SODEs over short time intervals, analogous to Taylor series for solutions of ODEs. Based on such expansions, the recent paper \cite{wanner2024higher} discusses several higher-order extensions of the standard first-order Kramers-Moyal formulas. The latter were originally proposed in \cite{boninsegna2018sparse} as a basis for extending SINDy to SODEs. Here we first recall from \cite{wanner2024higher} the main ingredients, presenting the estimation methods of interest in this work, to be later extended to SDDEs and then compared in the numerical Section \ref{sec:tests}. The extension to SDDEs is illustrated only at the end relying on an obvious argument. For a complete reference on stochastic-Taylor expansions for SODEs see \citep[Chapter 5]{kloeden1992stochastic}; as for SDDEs we refer to \citep{cao2015numerical,rosli2013systematicA,rosli2013stochastic}. We assume throughout that $f$ and $g$ in \eqref{eq:sode} (and then in \eqref{eq:sdde}) are sufficiently smooth with bounded derivatives to ensure the validity of the considered expansions as well as meaningful estimates of related remainders.

\bigskip
For the SODE \eqref{eq:sode} let
\begin{equation*}
L^{0} \coloneqq \sum_{i=1}^{n} f_{i}(x)\frac{\partial}{\partial x_{i}}
+\sum_{i,j=1}^{n} g_{i,j}(x)
\frac{\partial^{2}}{\partial x_{i} \partial x_{j}},\qquad
L^{j} \coloneqq \sum_{i=1}^{n} g_{i,j}(x) \frac{\partial}{\partial x_{i}},
\quad j=1,\dots,q,
\end{equation*}
be differential operators acting on sufficiently smooth maps $\mathbb{R}^{n}\to\mathbb{R}$, where subscripts denote components of vectors and matrices. $L^{0}$ captures both the drift and the second-order effects of diffusion (the It\^{o} correction), while $L^{j}$ represents the first-order effect of the $j$-th noise component.
For a specific test function $\phi$ the first few terms of the corresponding It\^{o}-Taylor expansion read
\begin{equation}\label{eq:itotaylor}
\setlength\arraycolsep{0.1em}\begin{array}{rcl}
\phi(X({t+\Delta t}))
&=&\displaystyle
\phi(X(t))
+L^{0}\phi(X(t))\Delta t
+\sum_{j=1}^{n} L^{j}\phi(X(t))\Delta W_{j}(t)\\
&&\displaystyle
+\frac{1}{2} (L^{0})^{2}\phi(X(t))\Delta t
+\sum_{j=1}^{n} L^{j}L^{0}\phi(X(t))I_{j,0}\\
&&\displaystyle
+\sum_{j=1}^{n} L^{0}L^{j}\phi(X(t))I_{0,j}
+\sum_{j,k=1}^{n} L^{j}L^k\phi(X(t))I_{j,k}
+\cdots
\end{array}
\end{equation}
where $\Delta W_{j}(t)\coloneqq W_{j}({t+\Delta t})-W_{j}(t)$ is the increment of the $j$-th component of the underlying Brownian motion and the terms $I$'s represent the multiple stochastic integrals
\begin{align*}
I_{j,0}&\coloneqq\int_t^{t+\Delta t}\int_t^{\sigma_{1}}\dd W_{j}({\sigma_{2}})\dd\sigma_{1},\\
I_{0,j}&\coloneqq\int_t^{t+\Delta t}\int_t^{\sigma_{1}}\dd\sigma_{2}\dd W_{j}({\sigma_{1}}),\\
I_{j,k}&\coloneqq\int_t^{t+\Delta t}\int_t^{\sigma_{1}}\dd W_{j}({\sigma_{2}})\dd W_{k}({\sigma_{1}}).
\end{align*}
\eqref{eq:itotaylor} can be used in weak form by taking the conditional expectation of $\phi(X({t+\Delta t}))$ given $X(t)$, or in strong form by working pathways with $\phi(X({t+\Delta t}))$ thus retaining stochasticity.
For instance, using the weak form with $\phi(x)\coloneqq x_{i}$ for some $i=1,\ldots,n$ leads to
\begin{equation*}
\mathbb{E}\left[x_{i}({t+\Delta t})\mid X(t)=x\right]
= x_{i} + f_{i}(x) \Delta t + O(\Delta t^{2})
\end{equation*}
from which the well-known Kramers-Moyal formula
\begin{equation}\label{eq:kmfw}
f_{i}(x)\approx\mathbb{E}\left[\frac{X_{i}({t+\Delta t})-X_{i}(t)}{\Delta t}\mid X(t)=x\right]
\end{equation}
follows as a component-wise estimate of the drift.
Similarly, for the covariance
\begin{equation}\label{eq:kmgw}
G_{i,j}(x)\approx\mathbb{E}\left[\frac{(X_{i}({t+\Delta t})-X_{i}(t))(X_{j}({t+\Delta t})-X_{j}(t))}{2\Delta t}\mid X(t)=x\right]
\end{equation}
follows by choosing $\phi(x)\coloneqq (x_{i}-X_{i}(t))(x_{j}-X_{j}(t))$. In the strong (pathways) version \eqref{eq:kmfw} and \eqref{eq:kmgw} become
\begin{equation}\label{eq:kmfs}
f_{i}(X(t))\approx\frac{X_{i}({t+\Delta t})-X_{i}(t)}{\Delta t}
\end{equation}
\begin{equation}\label{eq:kmgs}
G_{i,j}(X(t))\approx\frac{(X_{i}({t+\Delta t})-X_{i}(t))(X_{j}({t+\Delta t})-X_{j}(t))}{2\Delta t}.
\end{equation}
In either weak or strong form, in the sequel we denote by KM these standard first-order estimates \citep{kleinhans2012estimation,kleinhans2005iterative,risken1989methods}.

\bigskip
Following \cite{wanner2024higher}, improvements rely in general on higher-order finite-differences. In this work we consider only second-order estimates resorting respectively to forward, central and trapezoidal-like differences as illustrated next (only in the strong version for brevity).
\begin{itemize}
\item {\bf Forward Differences (FD).} A second-order forward finite-differences scheme follows for both drift and covariance by considering more terms in the It\^{o}-Taylor expansion \eqref{eq:itotaylor}:
\begin{equation}\label{eq:fdfs}
f_{i}(X(t))\approx\frac{4(X_{i}(t+\Delta t)-X_{i}(t))-(X_{i}(t+2\Delta t)-X_{i}(t))}{2\Delta t},
\end{equation}
\begin{equation}\label{eq:fdgs}
\setlength\arraycolsep{0.1em}\begin{array}{rcl}
G_{i,j}(X(t))&\approx&\displaystyle\frac{1}{4\Delta t}[4(X_{i}(t+\Delta t)-X_{i}(t))(X_{j}(t+\Delta t)-X_{j}(t))\\[2mm]
&&-(X_{i}(t+2\Delta t)-X_{i}(t))(X_{j}(t+2\Delta t)-X_{j}(t))].
\end{array}
\end{equation}
\item {\bf Central Differences (CD).} Similarly, second order can be reached in principle by
\begin{equation}\label{eq:cdfs}
f_{i}(X(t))\approx\frac{X_{i}(t+\Delta t)-X_{i}(t-\Delta t)}{2\Delta t},
\end{equation}
\begin{equation}\label{eq:cdgs}
\setlength\arraycolsep{0.1em}\begin{array}{rcl}
G_{i,j}(X(t))&\approx&\displaystyle\frac{(X_{i}(t+\Delta t)-X_{i}(t-\Delta t))(X_{j}(t+\Delta t)-X_{j}(t-\Delta t))}{4\Delta t}.
\end{array}
\end{equation}
Just ``in principle'' as the formulas above should lead to Stratonovich SODEs rather than It\^{o} ones \cite[Remark 2]{wanner2024higher}. We rather show in the experimental Section \ref{sec:tests} that they still work correctly, at least for the models considered therein.
\item {\bf Trapezoidal Rule (TR).} Finally, to avoid erroneous convergence of CD the authors of \cite{wanner2024higher} propose the following trapezoidal-like approximations:
\begin{equation}\label{eq:trfs}
\frac{f_{i}(X(t+\Delta t))+f_{i}(X(t))}{2}\approx\frac{X_{i}(t+\Delta t)-X_{i}(t)}{\Delta t},
\end{equation}
\begin{equation}\label{eq:trgs}
\setlength\arraycolsep{0.1em}\begin{array}{rcl}
G_{i,j}(X(t+\Delta t))+G_{i,j}(X(t))&\approx&\displaystyle
\frac{1}{2\Delta t}\{[2(X_{i}(t+\Delta t)-X_{i}(t))\\[2mm]
&&-\displaystyle
\Delta t(f_{i}(X(t+\Delta t))+f_{i}(X(t)))]\\[2mm]
&&\cdot\displaystyle
[2(X_{j}(t+\Delta t)-X_{j}(t))\\[2mm]
&&-\displaystyle
\Delta t(f_{j}(X(t+\Delta t))+f_{j}(X(t)))]\}.
\end{array}
\end{equation}
\end{itemize}
All the KM, FD, CD and TR estimates are employed and compared in Section \ref{sec:tests}. As discussed in Section \ref{sec:estimates}, resorting to their weak or strong form depends on the specific strategy to be adopted in view of the available data.

\bigskip
We close this section by observing that extension to SDDEs follows straightforwardly for the It\^{o}-Taylor expansions just by formally replacing the state variable $X(t)\in\mathbb{R}^{n}$ of the SODE \eqref{eq:sode} with an augmented state variable $Z(t)\coloneqq(X(t),X(t-\tau))\in\mathbb{R}^{2n}$ from the SDDE \eqref{eq:sdde}. As far as the related weak estimates of drift and covariance are concerned, we stress that now conditioning must include both current and delayed states for consistent parameter estimation due to the non-Markovian nature of SDDEs. As an instance \eqref{eq:kmfw} and \eqref{eq:kmgw} become
\begin{equation*}
f_{i}(z)\approx\mathbb{E}\left[\frac{X_{i}({t+\Delta t})-X_{i}(t)}{\Delta t}\mid X(t)=z_{1},\ X(t-\tau)=z_{2}\right],
\end{equation*}
\begin{equation*}
G_{i,j}(z)\approx\mathbb{E}\left[\frac{(X_{i}({t+\Delta t})-X_{i}(t))(X_{j}({t+\Delta t})-X_{j}(t))}{2\Delta t}\mid X(t)=z_{1},\ X(t-\tau)=z_{2}\right],
\end{equation*}
by choosing $\phi(z)=z_{1,i}$ in the It\^{o}-Taylor expansion. Above, the block partition of $z\in\mathbb{R}^{2n}$ is intended as $(z_{1},z_{2})\in\mathbb{R}^{n}\times\mathbb{R}^{n}$ with $z_{1}$ representing the current time value and $z_{2}$ the delayed counterpart.
\section{Sparse Identification of Nonlinear Dynamics with SINDy}\label{sec:sindy}
SINDy is a data-driven framework designed to uncover governing equations directly from observed time-series data. It leverages the assumption that most dynamical systems have an inherently sparse structure: their governing equations typically involve right-hand sides (RHSs) with a limited number of significant terms. This sparsity assumption allows the construction of accurate and interpretable models when the true underlying equations are unknown or when dimensionality reduction is required. We first summarize the original approach for ODEs in Section \ref{sec:sindy_ode}. Then we briefly illustrate the main aspects of the recent extensions to DDEs in Section \ref{sec:sindy_dde} and to SODEs in Section \ref{sec:sindy_sode}. Finally we discuss the approach followed in this work for tackling SDDEs in Section \ref{sec:sindy_sdde}.
\subsection{SINDy algorithm for ODEs}\label{sec:sindy_ode}
For general references see the pioneering work~\citep{bpk16} as well as~\cite{bk19,champ19}.

Consider a dynamical system governed by an ODE $x'=f(x)$ with unknown RHS
$f:\mathbb{R}^{n}\to\mathbb{R}^{n}$. Suppose to know or measure the state $x(t)$ at equidistant or even irregularly spaced discrete time points $t_{1},\dots,t_{m}$. SINDy approximates each component $f_{i}(x)$, $i = 1, \dots, n$, as a sparse linear combination of a predefined set of candidate basis functions $\{\theta_j(x)\}_{j=1}^{p}$:
\begin{equation}\label{eq:fi}
    f_{i}(x) = \sum_{j=1}^{p}\theta_j(x)\xi_{j,i}.
\end{equation}
The choice of basis functions $\theta_j(x)$ significantly affects the effectiveness of the approach (common basis functions include polynomials, trigonometric functions and other nonlinear functions that may reflect the anticipated system behaviour or prior knowledge).

To estimate the coefficients $\xi_{j,i}$, the collected data are organized as
\begin{equation*}
    \mathbf{X} \coloneqq \begin{bmatrix}
        x_{1}(t_{1}) &  \dots & x_n(t_{1}) \\      
        \vdots & \ddots  & \vdots \\ 
        x_{1}(t_{m})  & \dots & x_n(t_{m})
    \end{bmatrix}\in \mathbb{R}^{m \times n},\qquad 
    \mathbf{X'} \coloneqq \begin{bmatrix} 
    x_{1}'(t_{1})  & \dots & x_n'(t_{1}) \\
    \vdots & \ddots  & \vdots \\
    x_{1}'(t_{m})  & \dots & x_n'(t_{m})
    \end{bmatrix}\in \mathbb{R}^{m \times n},
\end{equation*} 
where the derivatives may either be directly measured or estimated numerically~\citep{chtd11}. A library matrix
\begin{equation*}
\Theta(\mathbf{X}) \coloneqq \begin{bmatrix}
1 & \mathbf{X} & \mathbf{X}^{2} & \cdots & \mathbf{X}^d & \cdots & \sin(\mathbf{X}) & \cos(\mathbf{X}) & \cdots
\end{bmatrix}\in \mathbb{R}^{m \times p}.
\end{equation*}
is then constructed by evaluating the candidate basis functions at all time points. Above, $\mathbf{X}^d$ denotes a matrix containing all possible polynomial combinations of degree $d$ in the components of $x$. 
The total number $p$ of basis functions depends on the chosen function types and their complexity (e.g., a polynomial library of degree $d$ include $p = \binom{n + d}{d}$ basis functions).

Determining the coefficients $\xi_{j,i}$ involves solving $\mathbf{X'} = \Theta(\mathbf{X})\Xi$ for $\Xi\coloneqq [\xi_{1}\; \xi_{2}\; \cdots\; \xi_n]\in\mathbb{R}^{p\times n}$ where each column vector $\xi_{i}\coloneqq(\xi_{1,i},\ldots,\xi_{p,i})\in\mathbb{R}^{p}$ corresponds to the sparse representation \eqref{eq:fi}.
Each $\xi_{i}$ is determined by independently solving\begin{equation}\label{eq:optknown}
\xi_{i} =  \argmin_{\xi \in \mathbb{R}^p}\left(\left\|\mathbf{X}_{i}'-\Theta(\mathbf{X})\xi\right\|_{2}+\lambda\left\|\xi\right\|_{1}\right) ,
\end{equation}
with sparse regression algorithms (such as STLS~\citep{bpk16} or LASSO  \cite{lasso96}). Above $\mathbf{X}_{i}'$ denotes the $i$-th column of $\mathbf{X'}$ and the hyperparameter $\lambda \geq 0$ controls the balance between fitting accuracy (the first term, in $\|\cdot\|_{2}$) and sparsity (the second term,  in $\|\cdot\|_{1}$).
\subsection{SINDy for DDEs}\label{sec:sindy_dde}
Recent extensions to DDEs such as $x'(t)=f(x(t),x(t-\tau))$ for a single constant delay $\tau>0$ have been proposed in \citep{bbt24,kopeczi23,pec24,sandoz2023sindy}. The following summary refers mainly to \citep{bbt24,pec24}. The core of SINDy for DDEs is the construction of an augmented library $\Theta(\mathbf{X},\mathbf{X}_{\tau})$ including delayed time samples $x(t_{i}-\tau)$, $i=1\ldots,m$. If these values are not at disposal they can be provided by piecewise linear interpolation of the available data $\mathbf{X}$. Then the SINDy framework extends by solving $\mathbf{X'} = \Theta(\mathbf{X},\mathbf{X}_{\tau})\Xi(\tau)$ via sparse regression, where now the sparse coefficients $\Xi(\tau)$ depend on the value of $\tau$, assumed to be known in principle. The case of multiple delays can be treated straightforwardly. Different is the realistic situation where also the delays are unknown (possibly both in number and values). The main idea is to minimize the reconstruction error $\tau\mapsto\|\mathbf{X'}-\Theta(\mathbf{X},\mathbf{X}_{\tau})\Xi(\tau)\|_{2}$ as  a function of the unknown delay(s) $\tau=(\tau_{1},\ldots,\tau_{k}=:\bar\tau)$ by using an ``external'' optimizer (not to be confused with the ``internal'' minimization of SINDy). In \cite{pec24} the use of Bayesian tools has largely improved the brute force approach originally followed in \cite{sandoz2023sindy}. In \cite{bbt24} a further substantial improvement has been reached by employing particle swarm techniques \citep{bonyadi2017,kennedy1995,shi1998}. Moreover, in \cite{bbt24} a different approach has been proposed, based on first discretizing the unknown DDE into a finte number of ODEs, to then apply the standard SINDy algorithm, possibly with external optimization. The great advantage of this ``pragmatic'' approach is that of having to deal only with the maximum delay $\bar\tau$, while the presence of possible multiple intermediate delays $(\tau_{1},\ldots,\tau_{k-1})$ is hidden in the collocation process that reduce the DDE to ODEs \cite{breda2016pseudospectral}. In this case the augmented library includes samples at time-shifted collocation nodes rather delayed samples. Both approaches can be extended by including stochastic terms as explained in Section \ref{sec:sindy_sdde}.
\subsection{SINDy for SODEs}\label{sec:sindy_sode}
The application of SINDy to SODEs such as \eqref{eq:sode} has been first considered in \cite{boninsegna2018sparse} and more recently in \cite{wanner2024higher} (see also the references therein and \citep{lenfesty2025uncovering,wang2022data}). The main difference from the use on deterministic differential equations, where data $\mathbf{X}'$ are either available or can be easily recovered from $\mathbf{X}$, is the need for reasonable data for drift and diffusion. In view of understanding just the principle underlying the extension of SINDy to SODEs, for the moment we can assume that such data are somehow reconstructed or estimated in a non-parametric fashion and are then available as $\mathbf{f}(\mathbf{X})\in\mathbf{R}^{m\times n}$ for drift and $\mathbf{G}(\mathbf{X})\in\mathbf{R}^{m\times n^{2}}$ for diffusion in the form of vectorized covariance. It is then enough to choose two possibly different libraries of candidate functions
$\Theta_{f}\in\mathbb R^{m\times p_{f}}$ and
$\Theta_{G}\in\mathbb R^{m\times p_{G}}$ and solving independently
\begin{align}
  \alpha_{i} 
  &= \argmin_{\alpha\in\mathbb R^p}
    \left\|\mathbf{f}_{i}(\mathbf{X}) - \Theta_{f}(\mathbf{X})\alpha \right\|_{2}
    + \lambda_{f}\|\alpha\|_{1},\quad i=1,\dots,n,\label{eq:sindy_sodef}\\
  \beta_{i,j}
  &= \argmin_{\beta\in\mathbb R^p}
    \left\|\mathbf{G}_{i,j}(\mathbf{X}) - \Theta_{G}(\mathbf{X})\beta\right\|_{2}
    + \lambda_{G}\|\beta\|_{1},\quad i,j=1,\dots,n.\label{eq:sindy_sodeg}
\end{align}
To note the freedom of choosing also different sparse promoting parameters among drif and diffusion. We demand to Section \ref{sec:estimates} a detailed discussion on how obtaining $\mathbf{f}(\mathbf{X})$ and $\mathbf{G}(\mathbf{X})$ following three different approaches. First we deal with SINDy for SDDEs as the driving novelty of the present work.
\subsection{SINDy for SDDEs}\label{sec:sindy_sdde}
To the best of the authors' knowledge the only work on SINDy for SDDEs is the recent paper \cite{han2024approximation}. The approach proposed therein heavily relies on assuming a single small delay in order to approximate SDDEs by SODEs via Taylor expansion. As such, the methodology is not suitable for generic delay values. Moreover, the stochastic integral is restricted to the Stratonovich type and an extension to multiple delays seems far from being trivial.

Here instead we pursue a general extension of SINDy to SDDEs by combining the arguments of both Section \ref{sec:sindy_dde} for DDEs and Section \ref{sec:sindy_sode} for SODEs. Indeed, the key observation is that in principle the presence of delay terms asks for an augmented library including delayed samples $\mathbf{X}_{\tau}$, while the presence of stochastic terms requires suitable strategies to recover data $\mathbf{f}(\mathbf{X})$ for drift and $\mathbf{G}(\mathbf{X})$ for diffusion through covariance. As such the two issues appear basically unrelated and this would lead to extend \eqref{eq:sindy_sodef} and \eqref{eq:sindy_sodeg} as
\begin{align*}
  \alpha_{i} 
  &= \argmin_{\alpha\in\mathbb R^p}
    \left\|\mathbf{f}_{i}(\mathbf{X}) - \Theta_{f}(\mathbf{X},\mathbf{X}_{\tau})\alpha \right\|_{2}
    + \lambda_{f}\|\alpha\|_{1},\quad i=1,\dots,n,\\
  \beta_{i,j}
  &= \argmin_{\beta\in\mathbb R^p}
    \left\|\mathbf{G}_{i,j}(\mathbf{X}) - \Theta_{G}(\mathbf{X},\mathbf{X}_{\tau})\beta\right\|_{2}
    + \lambda_{G}\|\beta\|_{1},\quad i,j=1,\dots,n,
\end{align*}
together with all the possible variations in view of multiple or unknown delays, as well as the pragmatic reduction of the delay dependence via collocation as mentioned in Section \ref{sec:sindy_dde}. Nevertheless, a more trustworthy approach should instead consider the explicit dependence on the delayed samples also for the reconstruction of data for drift and diffusion, hence
\begin{align}
  \alpha_{i} 
  &= \argmin_{\alpha\in\mathbb R^p}
    \left\|\mathbf{f}_{i}(\mathbf{X},\mathbf{X}_{\tau}) - \Theta_{f}(\mathbf{X},\mathbf{X}_{\tau})\alpha \right\|_{2}
    + \lambda_{f}\|\alpha\|_{1},\quad i=1,\dots,n,\label{eq:sindy_sddef}\\
  \beta_{i,j}
  &= \argmin_{\beta\in\mathbb R^p}
    \left\|\mathbf{G}_{i,j}(\mathbf{X},\mathbf{X}_{\tau}) - \Theta_{G}(\mathbf{X},\mathbf{X}_{\tau})\beta\right\|_{2}
    + \lambda_{G}\|\beta\|_{1},\quad i,j=1,\dots,n.\label{eq:sindy_sddeg}
\end{align}
Indeed, as discussed at the end of Section \ref{sec:ito_taylor} in view of It\^{o}-Taylor expansions for SDDEs, in general weak forms require conditioning also w.r.t. past values and thus \eqref{eq:sindy_sodef} and \eqref{eq:sindy_sodeg} fit better into the framework. The issue becomes evident in the forthcoming Section \ref{sec:estimates} about one of the possible strategies to recover the necessary data.
\section{Estimation of drift and diffusion from data}\label{sec:estimates}
In this section we deal with obtaining the data $\mathbf{f}(\mathbf{X})$ and $\mathbf{G}(\mathbf{X})$ as just anticipated in Section \ref{sec:sindy_sdde}. In Section \ref{sec:ito_taylor} general moment-based estimators are described relying on It\^{o}-Taylor expansions. Independently of the specific choice or order, KM, FD, CD and TR methods all face a fundamental challenge: how to get sufficiently-many data at each time instant to ensure reliable estimates of drift and diffusion. The answer depends, among other factors, on whether temporal data are available for a single trajectory alone, or samples from multiple independent realizations are at disposal (realistically from repeated experiments). In the first case the approach we follow, hereafter named ``approach A'' or simply ``A'' for brevity, relies on creating more realizations via suitable tools and is described in Section \ref{sec:A}. In the second case, ``approach B'' or simply ``B'' in the sequel, we discuss two different strategies, namely B1 and B2, described respectively in Section \ref{sec:B1} and in Section \ref{sec:B2}. The key difference between B1 and B2 is whether averages are taken before or after the application of SINDy's sparse regression.
Then in Section \ref{sec:tests} we perform a thorough experimental comparison on different models of SDDEs, illustrating different trade-off's between reconstruction error, computational cost, data requirements and implementation complexity.

First to proceed we underline that most of the adopted techniques we refer to in the sequel are taken from the recent literature but none of them deal with the presence of delays. We thus extend their application by considering when necessary an augmented state $Z(t)\coloneqq (X(t),X(t-\tau))\in\mathbb{R}^{2n}$ as already discussed at the end of Section \ref{sec:ito_taylor}. Nevertheless, we are aware that the general claims about the properties and effectiveness of such tools that hold in the case without delays are not rigorously proved for SDDEs in this work (as this would lead us out of scope). We reserve to dive into this direction in the next future, relying for the moment on the experimental results of Section \ref{sec:tests}.
\subsection{Approach A: more realizations from conditional expectations}\label{sec:A}
In this section we resort to the methodology discussed recently in \cite{chen2022non}. The underlying basis is used also in the first work \cite{boninsegna2018sparse} on SINDy for SODEs by employing KM estimates.

Let $\delta t$ be the stepsize chosen to represent the Brownian motion $W$ underlying \eqref{eq:sdde} and assume that samples of a single path $X(t)$ are available with time stepsize $\Delta t=K\delta t$ for some integer $K\geq1$ (for our simulation we have used $K=1$ in all the results presented in Section \ref{sec:tests}). If not, one may think at interpolated samples, also to get delayed samples and hence $Z(t)$. Given the relevant time instants $t_{i}$, $i=1,\ldots,m$ and a positive integer $M$, for each $i=1,\ldots,m-1$ the proposed methodology first requires the following three steps:
\begin{enumerate}
\item[A(i)] initialize $M$ independent sample paths $X^{(k)}(t)$, $k=1,\ldots,M$, of \eqref{eq:sdde} each starting at time $t_{i}$  from $Z(t_{i})$;
\item[A(ii)] approximate a single step of each path at time $t_{i}+\Delta t$ by using a numerical time-integrator, e.g., Euler–Maruyama as adapted to \eqref{eq:sdde} \cite{buckwar2000introduction}, viz.
\begin{equation*}
X^{(k)}_{i+1}=X^{(k)}_{i}+f(X^{(k)}_{i},X^{(k)}_{i-\kappa})\Delta t+g(X^{(k)}_{i},X^{(k)}_{i-\kappa})(W^{(k)}(t_{i+1})-W^{(k)}(t_{i})),
\end{equation*}
where we implicitly assume that $\tau$ is an integer multiple $\kappa$ of $\Delta t$;
\item[A(iii)] record the increment $\Delta X^{(k)}_{i}\coloneqq X^{(k)}_{i+1}-X^{(k)}_{i}$.
\end{enumerate}
Then estimates of drift and diffusion (again through covariance) are obtained as
\begin{equation*}
      \mathbf{f}(Z({t_{i}}))
  = \frac{1}{M\Delta t}
    \sum_{k=1}^{M} \Delta X_{i}^{(k)},
\end{equation*}
\begin{equation*}
  \mathbf{G}(Z({t_{i}}))
  = \frac{1}{\Delta t}
    \left(\frac{1}{M}\sum_{k=1}^{M}
           \Delta X_{i}^{(k)}\left(\Delta X_{i}^{(k)}\right)^{T}
         - \mathbf{f}(Z({t_{i}}))\mathbf{f}(Z({t_{i}}))^{T} (\Delta t)^2\right).
\end{equation*}
As $M\to\infty$ and $\Delta t\to0$, these estimates are expected to converge to the true values $f(Z(t_{i}))$ and $G(Z(t_{i}))$ (which is indeed the case without delay \cite{chen2022non}). In practice one chooses $K$ large enough ($\delta t\ll\Delta t$) to control discretization bias and $M$ large enough to reduce sampling variance.  Higher‐order schemes (as, e.g., Milstein \cite{mil1975approximate} and stochastic Runge–Kutta \cite{kloeden1992stochastic}) or finite‐time corrections \cite{ragwitz2001indispensable} can be used when $\Delta t$ is not negligible.

A final comment is fundamental, also in view of approaches B1 and B2 as described next: the above procedure requires an assumed specification or initial approximation of both $f$ and $g$, otherwise step A(ii) is not possible. In practice, one supplies these functions via a suitable parametric family, a non-parametric initializer or prior structural information and then refines them iteratively until a desired tolerance \cite{chen2022non}. In Section \ref{sec:tests} we adopt the true forms of $f$ and $g$ for the sake of avoiding complexities unrelated to the current scope. 
\subsection{Approach B1: multiple trajectories with pre-regression averaging}
\label{sec:B1}
An alternative to approach A can be followed when a large dataset from several independent realizations of the system is available.

Consider $M$ independent trajectories of \eqref{eq:sdde}, each sampled at $t_{i}$, $i=1,\ldots,m$, with uniform time step $\Delta t$. Let $Z^{(k)}({t_{i}})$ be the augmented state of the $k$-th trajectory at time $t_{i}$ and $\Delta X^{(k)}({t_{i}}) = X^{(k)}({t_{i+1}}) - X^{(k)}({t_{i}})$ be the relevant increment of the current time state. To estimate drift and diffusion at a specific point $z\in \mathbb{R}^{2n}$ in the augmented state space we aggregate data points from all trajectories falling within a predefined neighborhood of $z$ \citep{bandi2003fully,mohammadi2024nonparametric}. Given a tolerance $\varepsilon>0$ let $N_{\varepsilon}(z)$ denote the cardinality of the set $S_{\varepsilon}(z)\coloneqq\{(k,i)\ :\ \|Z^{(k)}(t_{i}) - z\| \leq \varepsilon\}$ of such observations. 
This effectively performs a kernel-like local averaging, where 
$\varepsilon $ controls the trade-off between statistical variance (too few points if $\varepsilon$ is small) and bias (over-smoothing if $\varepsilon$ is large). This and similar techniques are widely used in non-parametric drift and diffusion estimation methods for stochastic processes for instance, through kernel density estimation or binning schemes in continuous-time models. The required estimates are finally obtained by averaging the observed increments and quadratic variations, respectively, over all data points in the neighborhood of $z$, viz.
\begin{align}
\mathbf{f}(z)
&= \frac{1}{N_{\varepsilon}(z)\Delta t}
\sum_{k=1}^{M}
\sum_{i:(k,i)\in S_{\varepsilon}(z)}
\Delta X_{t_{i}}^{(k)},
\label{eq:mt-sdde-drift-revised}\\[6pt]
\mathbf{G}(z)
 &= \frac{1}{N_{\varepsilon}(z)\Delta t}
\sum_{k=1}^{M}
\sum_{i:(k,i)\in S_{\varepsilon}(z)}
\Delta X_{t_{i}}^{(k)}
(\Delta X_{t_{i}}^{(k)})^{T}.
\label{eq:mt-sdde-diffusion-revised}
\end{align}
The estimates presented as such are based on KM in an ensemble setting. Generalization to FD, CD and TR (presented in Section \ref{sec:ito_taylor}) is straightforward. Convergence to $f(z)$ and $G(z)$ is expected as $M \to \infty$ (or as the total number of observations $N_{\varepsilon}(z)$ increases) also for SDDEs. The accuracy depends on balancing the bias (from finite $\Delta t$ and neighborhood size $\varepsilon$) and variance (from limited amount of data).
As a final comment we observe that this approach facilitates a sparse, data-driven reconstruction of the underlying non-Markovian stochastic dynamics of SDDEs without requiring explicit knowledge of the functional forms of $f$ and $g$ beforehand. Moreover, we remark once more that B1 is characterized by taking averages before the application of SINDy's sparse regression, contrary to B2 as discussed next.
\subsection{Approach B2: multiple trajectories with post-regression averaging }
\label{sec:B2}
A second strategy for leveraging an ensemble of $M$ trajectories consists in applying SINDy's sparse regression to each trajectory individually and then averaging the identified model parameters. This method, described in the experimental setup of \cite{wanner2024higher}, contrasts with B1 by averaging in the parameter space after regression, rather than in the data space before regression. Briefly, if $\alpha_{i}^{(k)}$ and $\beta_{i,j}^{(k)}$ are the sparse coefficients identified by SINDy from the $k$-th trajectory through respectively \eqref{eq:sindy_sodef} and \eqref{eq:sindy_sodeg}, the final model coefficients are retrieved by averaging across the entire ensemble:
\begin{equation*}
\alpha_{i}=\frac{1}{M}\sum_{k=1}^{M}\alpha_{i}^{(k)},\qquad\beta_{i,j}= \frac{1}{M} \sum_{k=1}^{M} \beta_{i,j}^{(k)}.
\end{equation*}
Practically B2 creates an ensemble of models, with the final model being the average of this ensemble. The underlying assumption is that while each individually identified model may contain variance due to the specific realization of the stochastic process, averaging over many such models yields a result that is expected to converge to the true underlying system dynamics.
\section{Numerical experiments and results}\label{sec:tests}
In this section we perform a thorough computational analysis by collecting several experiments on the three SDDEs described below. The experiments consist in applying the approaches A, B1 and B2 discussed in Section \ref{sec:estimates}, each of them using all of KM, FD, CD and TR as underlying estimators from Section \ref{sec:ito_taylor}. We recall that $M$ denotes both the numbers of paths realized with A through Euler-Maruyama simulations, as well as the cardinality of the ensemble of available trajectories for B1 and B2. A full comparison is first performed for all the three examples by assuming $M=1\,000$ and $\varepsilon=10^{-4}$ for B1. Further tests regard limited comparisons to separately assess the role of $M$ and $\varepsilon$ also w.r.t. CPU time and errors. As for the latter, we consider both the absolute errors $|\bar\xi_{i}-\xi_{i}|$ on the $i$-th sparse coefficient of the chosen SINDy library, $i=1,\ldots,p$, as well as the RMSE on relevant trajectories (they regard both trajectories for the state $X(t)$ as well as for separated drift and diffusion dynamics reconstruction). Finally, data are artificially provided by numerical simulation via Euler-Maruyama and we use STLS~\citep{bpk16} for sparse regression because of its robustness compared to other methods. All tests are performed on a Windows 11 system (CPU 5 GHz, RAM 16 GB) using MATLAB implementations for SINDy (MATLAB version R2024b). Before discussing the obtained results in Section \ref{sec:discussion} we present the three specific SDDEs that we use as models to test.

\bigskip
As a first model we consider from \citep{babasola2023stochastic,wu2012stochastic} the stochastic version with multiplicative noise
\begin{equation}\label{eq:stoch_delay_logistic}
    \dd X(t) = \alpha X(t)\bigl(1 - X(t-\tau)\bigr)\dd t + \sigma X(t) \dd W(t)
\end{equation}
of the so-called delay logistic or Hutchinson equation \cite{Hutchinson1948}. \eqref{eq:stoch_delay_logistic} models the evolution of a single ($n=1$) population $X(t)$ at time $t$ with intrinsic growth rate $\alpha=2$, time delay $\tau=1$ and noise intensity $\sigma=0.4$ with underlying scalar ($q=1$) Brownian motion. Data are obtained by simulating over the time window $[0,20]$ with time step $\Delta t=0.01$ via Euler–Maruyama, prescribing the history $X(s)=\cos(s)$, $s\in[-\tau,0]$. 80\% of the samples are used for training and 20\% for validation. SINDy is applied with polynomial libraries of degree $d=2$ for both drift and diffusion and $\lambda=0.025$ is used. Results of the full comparison are collected in Table \ref{tab:comparison1} and discussed in Section \ref{sec:discussion}.
\begin{sidewaystable}[htbp]
\centering
\scriptsize
\begin{tabular}{ll c ccc ccc ccc ccc}
\toprule
& library & true & \multicolumn{12}{c}{$|\xi_{i}-{\bar{\xi_{i}}}|$} \\
\cmidrule(r){4-15} 
& terms & coeff. & \multicolumn{3}{c}{KM} & \multicolumn{3}{c}{FD} & \multicolumn{3}{c}{CD} & \multicolumn{3}{c}{TR} \\
\cmidrule(r){4-6} \cmidrule(r){7-9} \cmidrule(r){10-12} \cmidrule(r){13-15} 
& &  & $A$  & B1& B2 & $A$ &  B1 & B2 & $A$ &  B1 & B2 & $A$  & B1 & B2 \\
\midrule
\multirow{2}{*}{drift} 
& $X(t)$ & $\xi_{2}=2$ & {$2.13e^{-3}$} & \boldmath{$1.29e^{-4}$} & {$5.76e^{-1}$} & $2.64e^{-2}$  & \boldmath{$1.28e^{-4}$} & {$5.54e^{-1}$} & {$1.11e^{-2}$} &  \boldmath{$3.87e^{-4}$} & {$3.36e^{-2}$} &\boldmath{$2.13e^{-3}$}  & \boldmath{$2.13e^{-3}$} & {$1.02e^{-1}$}\\
& $X(t)X(t-\tau)$ & $\xi_{5}=2$ & {$5.07e^{-3}$} & \boldmath{$6.68e^{-5}$} & {$7.47e^{-2}$} & {$2.53e^{-2}$}  & \boldmath{$7.14e^{-5}$} & {$6.79e^{-2}$} & {$1.05e^{-2}$}  & \boldmath{$2.08e^{-4}$} & {$6.57e^{-2}$} &{$2.10e^{0}$}  & \boldmath{$1.17e^{-3}$}& {$1.60e^{-2}$} \\
\midrule
\multicolumn{15}{c}{\textbf{same diffusion}} \\
\midrule
diffusion
& $X^{2}(t)$ & $\xi_{4}=0.16$ & {$6.39e^{-3}$}  & \boldmath{$7.54e^{-9}$} & {$9.17e^{-3}$} & {$6.39e^{-3}$}  & \boldmath{$7.54e^{-9}$}  & {$9.17e^{-3}$} & {$6.39e^{-3}$} &  \boldmath{$7.54e^{-9}$}   & {$9.17e^{-3}$} & {$6.39e^{-3}$}   & \boldmath{$7.54e^{-9}$}  & {$9.17e^{-3}$} \\

\midrule
\multicolumn{15}{c}{\textbf{performance metrics}} \\
\midrule
CPU [s] & \multicolumn{2}{c}{} & $\approx 4$  & $\approx 10$ & $\approx 10$ & $\approx 4$  & $\approx 10$ &$\approx 10$   &  $\approx 4$  & $\approx 10$  & $\approx 10$ & $\approx 4$   & $\approx 10$ & $\approx 10$ \\
\midrule
& \multicolumn{2}{c}{full} & {$2.99e^{-1}$}  & \boldmath{$1.10e^{-3}$} & {$8.79e^{-1}$} & {$3.43e^{-1}$}  & \boldmath{$4.29e^{-4}$} & $9.80e^{-1}$ & {$3.31e^{-1}$}  & \boldmath{$1.30e^{-3}$} & {$1.10e^{0}$} &{$3.73e^{0}$}  & \boldmath{$1.43e^{-3}$}& {$1.68e^{0}$} \\
RMSE & \multicolumn{2}{c}{drift} & {$1.11e^{-1}$}  & \boldmath{$ 1.13e^{-4}$} & {$3.79e^{-1}$} & $2.18e^{-2}$  & \boldmath{$9.75e^{-5}$}  & $4.68e^{-1}$ & {$9.00e^{-3}$}  & \boldmath{$3.07e^{-4}$}  & {$1.37e^{-1}$} & {$5.47e^{0}$} & \boldmath{$1.70e^{-3}$} & {$1.13e^{0}$} \\
& \multicolumn{2}{c}{diffusion}  & $9.99e^{-2}$ & \boldmath{$4.92e^{-6}$} &  {$2.43e^{-1}$} & {$9.97e^{-2}$}  & \boldmath{$2.85e^{-6}$} &  $2.21e^{-1}$ & {$9.97e^{-2}$}  & \boldmath{$1.53e^{-5}$} & {$1.97e^{-1}$} &{$9.99e^{-2}$}  & \boldmath{$1.22e^{-5}$}& {$2.20e^{-1}$}  \\

\midrule
\multicolumn{15}{c}{\textbf{different diffusion}} \\
\midrule
diffusion
& $X^{2}(t)$ & $\xi_{4}=0.16$ & {$6.39e^{-3}$}  & \boldmath{$7.54e^{-9}$} & {$8.31e^{-2}$} & {$7.91e^{-2}$}  & \boldmath{$1.50e^{-5}$} &  {$7.92e^{-2}$}& $8.48e^{-2}$ & \boldmath{$7.76e^{-3}$} &  {$9.47e^{-2}$} &$1.84e^{-1}$    &\boldmath{$2.90e^{-6}$} & $6.95e^{-2}$  \\

\midrule
\multicolumn{15}{c}{\textbf{performance metrics}} \\
\midrule
CPU [s] & \multicolumn{2}{c}{} & $\approx 4$ &  $\approx 10$ & $\approx 10$  & $\approx 4$  & $\approx 10$ & $\approx 10$ &  $\approx 4$  & $\approx 10$  &$\approx 10$  & $\approx 4$ &  $\approx 10$ & $\approx 10$ \\
\midrule
& \multicolumn{2}{c}{full} & {$2.99e^{-1}$}  & \boldmath{$1.10e^{-3}$} & $9.43e^{-1}$ & {$5.25e^{-1}$}  &  \boldmath{$6.00e^{-3}$} & {$1.08e^{0}$}&  $6.49e^{-1}$  & \boldmath{$2.99e^{-2}$}  &  {$1.89e^{0}$} & $2.62e^{0}$   & \boldmath{$3.00e^{-2}$}& {$9.83e^{-1}$}\\
RMSE & \multicolumn{2}{c}{drift} & {$1.11e^{-2}$}  & \boldmath{$ 1.13e^{-4}$} & $1.15e^{0}$  & {$1.14e^{-2}$}  & \boldmath{$ 9.06e^{-4}$}  & {$3.59e^{-1}$} & \boldmath{$9.00e^{-3}$}  & {$2.32e^{-2}$} &{$1.56e^{-1}$} & {$5.47e^{0}$}  & \boldmath{$1.10e^{-3}$} &{$2.90e^{-1}$}  \\
& \multicolumn{2}{c}{diffusion} & $9.99e^{-2}$  & \boldmath{$4.92e^{-6}$} & $1.57e^{0}$ &  {$3.08e^{-1}$}  & \boldmath{$3.89e^{-5}$} & {$1.05e^{0}$} & {$2.69e^{-1}$}  & \boldmath{$1.54e^{-2}$} & {$8.97e^{-1}$} & {$6.21e^{-1}$} & \boldmath{$7.64e^{-6}$}& {$2.41e^{-1}$} \\
\bottomrule
\end{tabular}
\caption{comparison of errors on estimated coefficients, RMSE on trajectories and CPU time for the logistic SDDE \eqref{eq:stoch_delay_logistic} for both drift and diffusion, all estimation methods KM, FD, CD and TR from Section \eqref{sec:ito_taylor} and all data reconstruction strategies A, B1 and B2 from Section \ref{sec:estimates}. Diffusion estimation might remain the same (KM in particular, first part above) or change (all of KM, FD, CD and TR, second part below).}
\label{tab:comparison1}
\end{sidewaystable}
\begin{figure}
    \centering
    \includegraphics[width=1\linewidth]{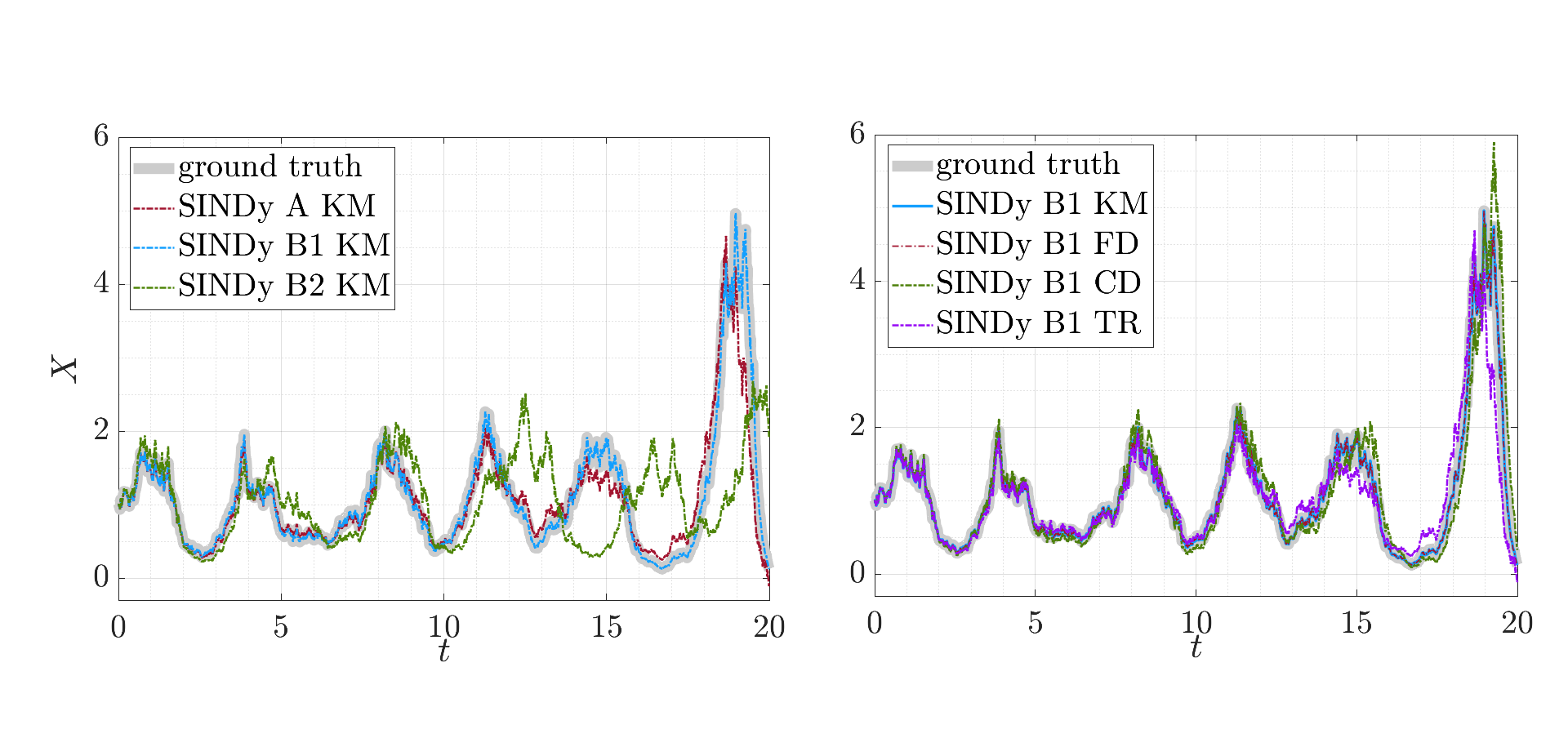}
    \caption{comparison of trajectories for the logistic SDDE \eqref{eq:stoch_delay_logistic}: ground truth and approaches A, B1 and B2 all with KM (left). Ground truth versus estimation methods KM, FD, CD and TR all with approach B1 (right). }\label{fig:log_traj}
\end{figure}

\bigskip
As a second model we consider from \citep{du2015stochastic}
the predator–prey interaction with delay and multiplicative noise
\begin{equation}\label{eq:pred_prey}
\left\{\setlength\arraycolsep{0.1em}\begin{array}{rcl}
\dd X_{1}(t)&=& X_{1}(t)[\alpha - \beta X_{1}(t) - \gamma X_{2}(t-\tau)]\dd t + \sigma_{1}X_{1}(t)\dd W_{1}(t)\\[2mm]
\dd X_{2}(t) &=& X_{2}(t)[-\delta + \kappa X_{1}(t-\tau)]\dd t + \sigma_{2}X_{2}(t)\dd W_{2}(t),
\end{array}\right.
\end{equation}
where $X_{1}(t)$ and $X_{2}(t)$ denote prey and predator populations, respectively, with parameters $\alpha=1$, $\delta=0.5$, $\beta=\gamma=\kappa=0.1$, $\sigma_{1}=\sigma_{2}=0.4$, $\tau=1$ and $W_{1},W_{2}$ independent scalar Wiener processes (here $n=q=2$). Again, data are obtained by simulating over the time window $[0,20]$ with time step $\Delta t=0.01$ via Euler–Maruyama, prescribing the constant history $X(s)\equiv(5,2)^{T}$. 80\% of the samples are used for training and 20\% for validation. SINDy is applied with polynomial libraries of degree $d=2$ for both drift and diffusion and $\lambda=0.025$ is used. Results of the full comparison are collected in Table \ref{tab:comparison2} and discussed in Section \ref{sec:discussion}.
\begin{sidewaystable}[htbp]
\centering
\scriptsize
\begin{tabular}{ll c ccc ccc ccc ccc}
\toprule
& library & true & \multicolumn{12}{c}{$|\xi_{i}-{\bar{\xi_{i}}}|$} \\
\cmidrule(r){4-15} 
& terms & coeff. & \multicolumn{3}{c}{KM} & \multicolumn{3}{c}{FD} & \multicolumn{3}{c}{CD} & \multicolumn{3}{c}{TR} \\
\cmidrule(r){4-6} \cmidrule(r){7-9} \cmidrule(r){10-12} \cmidrule(r){13-15} 
& &  & $A$  & B1 & B2 & $A$  & B1 & B2 & $A$  & B1  & B2 &$A$  & B1 & B2 \\
\midrule
 
\multirow{5}{*}{drift} 
& $X(t)$ & $\xi_{2}=1$ & {$1.28e^{-1}$}  & \boldmath{$9.20e^{-5}$} & {$2.27e^{0}$} &  {$6.90e^{-2}$}  & \boldmath{$5.19e^{-4}$} & $2.26e^{0}$ & {$1.16e^{-2}$} &  \boldmath{$5.52e^{-4}$} & $1.80e^{-1}$ & $3.23e^{-1}$ & \boldmath{$1.85e^{-5}$} & $4.51e^{0}$  \\
& $Y(t)$ & $\xi_{3}=0.5$ & {$1.42e^{-3}$}  & \boldmath{$2.86e^{-7}$} & {$1.13e^{0}$} & {$7.14e^{-2}$}  & \boldmath{$2.32e^{-4}$} & $1.55e^{0}$ & $1.78e^{-2}$ &  \boldmath{$2.56e^{-4}$} & $3.60e^{-1}$ & $7.86e^{-3}$ &  \boldmath{$1.38e^{-6}$} & $15.58e^{0}$ \\
& $X^{2}(t)$ & $\xi_{6}=0.1$ & {$6.16e^{-3}$}  & \boldmath{$ 8.38e^{-6}$} & {$1.65e^{-1}$} & {$6.59e^{-3}$}  & \boldmath{$3.98e^{-5}$} & $1.73e^{-1}$ & {$1.64e^{-3}$}  & \boldmath{$5.09e^{-5}$} & $7.88e^{-2}$  & {$1.68e^{-2}$} & \boldmath{$5.06e^{-5}$} & $6.33e^{-1}$ \\
& $X(t)Y(t-\tau)$ & $\xi_{9}=0.1$ & {$1.05e^{-2}$}  & \boldmath{$3.07e^{-6}$} & {$1.42e^{-1}$} & {$6.55e^{-3}$}  & \boldmath{$5.86e^{-5}$} & $1.45e^{-1}$ & {$7.63e^{-4}$} & \boldmath{$1.60e^{-5}$} & $2.20e^{-2}$ & $1.36e^{-1}$ & \boldmath{$3.79e^{-5}$} &  $5.60e^{-1}$  \\
& $Y(t)X(t-\tau)$ & $\xi_{11}=0.1$ & {$7.35e^{-4}$}  & \boldmath{$1.84e^{-8}$ } & {$5.07e^{-2}$} & {$1.08e^{-2}$}  & \boldmath{$2.35e^{-5}$} & $6.35e^{-2}$  & $3.07e^{-3}$  & \boldmath{$2.77e^{-5}$} & $1.15e^{-2}$  & $1.30e^{-2}$ & \boldmath{$1.45e^{-7}$} & $2.46e^{-1}$ \\
\midrule
\multicolumn{15}{c}{\textbf{same diffusion}} \\
\midrule
\multirow{2}{*}{ diffusion}
& $X^{2}(t)$ & $\xi_{6}=0.16$ & {$6.36e^{-4}$}  & \boldmath{$2.34e^{-6}$} & {$1.10e^{-2}$} & {$6.36e^{-4}$}  & \boldmath{$2.34e^{-6}$}  & {$1.10e^{-2}$}  & {$6.36e^{-4}$} & \boldmath{$2.34e^{-6}$}  & {$1.10e^{-2}$}  & {$6.36e^{-4}$} & \boldmath{$2.34e^{-6}$} & {$1.10e^{-2}$} \\
& $Y^{2}(t)$ & $\xi_{10}=0.16$ & $1.16e^{-4}$  & \boldmath{$1.06e^{-9}$} & {$3.55e^{-3}$} & $1.16e^{-4}$  & \boldmath{$1.06e^{-9}$}  & {$3.55e^{-3}$} & $1.16e^{-4}$  & \boldmath{$1.06e^{-9}$}  & {$3.55e^{-3}$}  & $1.16e^{-4}$ & \boldmath{$1.06e^{-9}$}  &  {$3.55e^{-3}$} \\
\midrule
\multicolumn{15}{c}{{\bf performance metrics}} \\
\midrule
CPU [s] & \multicolumn{2}{c}{} & $\approx 5$  & $\approx 30$ & $\approx 30$ & $\approx 5$  & $\approx 30$ & $\approx 30$ & $\approx 5$  & $\approx 30$ & $\approx 30$  & $\approx 5$& $\approx 30$ & $\approx 30$ \\
\midrule
& \multicolumn{2}{c}{full} & {$1.07e^{-1}$} &  \boldmath{$1.00e^{-3}$} & {$2.89e^{0}$} & {$4.84e^{-1}$} &  \boldmath{$4.80e^{-3}$} & {$2.13e^{0}$} & $1.50e^{-1}$ &  \boldmath{$7.20e^{-3}$} & $4.30e^{0}$ &   $1.70e^{0}$  & \boldmath{$5.37e^{-2}$} & $5.10e^{0}$ \\
RMSE & \multicolumn{2}{c}{drift} & {$3.10e^{-2}$}  & \boldmath{$ 8.80e^{-5}$} & {$2.66e^{0}$} & {$1.00e^{-1}$}  &  \boldmath{$5.40e^{-4}$} & {$8.89e^{0}$} & {$2.63e^{-2}$}  & \boldmath{$7.50e^{-4}$} & $1.27e^{0}$ &  $1.39e^{0}$  & \boldmath{$4.66e^{-4}$} & $2.03e^{1}$ \\
& \multicolumn{2}{c}{diffusion} & {$2.06e^{-2}$}  & \boldmath{$7.57e^{-6}$} & {$1.59e^{-1}$} &  {$2.28e^{-2}$}  & \boldmath{$6.62e^{-6}$} & {$8.87e^{-2}$} & {$2.28e^{-2}$} & \boldmath{$1.11e^{-5}$} & $6.93e^{-2}$ &  {$2.60e^{-2}$}   & \boldmath{$6.68e^{-5}$} &  $1.41e^{-1}$\\
\midrule
\multicolumn{15}{c}{\textbf{different diffusion}} \\
\midrule
\multirow{2}{*}{ diffusion}
& $X^{2}(t)$ & $\xi_{6}=0.16$ & {$6.36e^{-4}$} & \boldmath{$2.34e^{-6}$} & {$1.10e^{-2}$} & {$7.91e^{-2}$}  & \boldmath{$9.35e^{-6}$} & $7.13e^{-2}$ & $7.96e^{-2}$  & \boldmath{$3.13e^{-3}$} & $1.04e^{-1}$  &{$1.73e^{-3}$} & \boldmath{$1.30e^{-8}$}& {$1.10e^{0}$}\\
& $Y^{2}(t)$ & $\xi_{10}=0.16$ & $1.16e^{-4}$   & \boldmath{$1.06e^{-9}$} & {$3.55e^{-3}$} & $8.06e^{-2}$ & \boldmath{$2.32e^{-8}$} & $7.76e^{-2}$ & $7.97e^{-2}$  & \boldmath{$2.65e^{-4}$} & $1.05e^{-1}$  & {$7.61e^{-5}$} & \boldmath{$4.49e^{-9}$}& {$7.24e^{-1}$}\\
\midrule
\multicolumn{15}{c}{{\bf performance metrics}} \\
\midrule
CPU [s] & \multicolumn{2}{c}{} & $\approx 5$  & $\approx 30$ & $\approx 30$ & $\approx 5$  & $\approx 30$ & $\approx 30$ & $\approx 5$  & $\approx 30$ & $\approx 30$  & $\approx 5$& $\approx 30$ & $\approx 30$ \\
\midrule
& \multicolumn{2}{c}{full} & {$1.07e^{-1}$}  &  \boldmath{$1.00e^{-3}$} & {$2.89e^{0}$} & {$1.17e^{0}$}  &  \boldmath{$2.70e^{-3}$} & $2.21e^{0}$ & $1.27e^{0}$ &  \boldmath{$1.70e^{-2}$} & $1.14e^{0}$ & {$1.70e^{0}$} &  \boldmath{$5.80e^{-2}$} & $6.58e^{1}$ \\
RMSE & \multicolumn{2}{c}{drift} & {$3.10e^{-2}$}  & \boldmath{$ 8.80e^{-5}$} & {$2.66e^{0}$}  & {$1.00e^{-1}$}  & \boldmath{$5.76e^{-4}$} & $8.71e^{0}$ & {$2.63e^{-2}$} &  \boldmath{$2.34e^{-2}$} & $1.16e^{0}$ & {$1.39e^{0}$} & \boldmath{$3.90e^{-4}$} & $6.53e^{1}$\\
& \multicolumn{2}{c}{diffusion} & {$2.06e^{-2}$}  & \boldmath{$7.57e^{-6}$} & {$1.59e^{-1}$} & {$4.24e^{0}$}  & \boldmath{$3.38e^{-4}$} & $7.61e^{-1}$ & {$4.22e^{0}$} &  \boldmath{$6.06e^{-2}$} & {$6.72e^{-1}$} & {$5.42e^{-2}$} &  \boldmath{$1.55e^{-6}$} & $1.87e^{0}$ \\
\bottomrule
\end{tabular}
\caption{comparison of errors on estimated coefficients, RMSE on trajectories and CPU time for the predator-prey SDDE \eqref{eq:pred_prey} for both drift and diffusion, all estimation methods KM, FD, CD and TR from Section \eqref{sec:ito_taylor} and all data reconstruction strategies A, B1 and B2 from Section \ref{sec:estimates}. Diffusion estimation might remain the same (KM in particular, first part above) or change (all of KM, FD, CD and TR, second part below).}
\label{tab:comparison2}
\end{sidewaystable}
\begin{figure}
    \centering
    \includegraphics[width=1\linewidth]{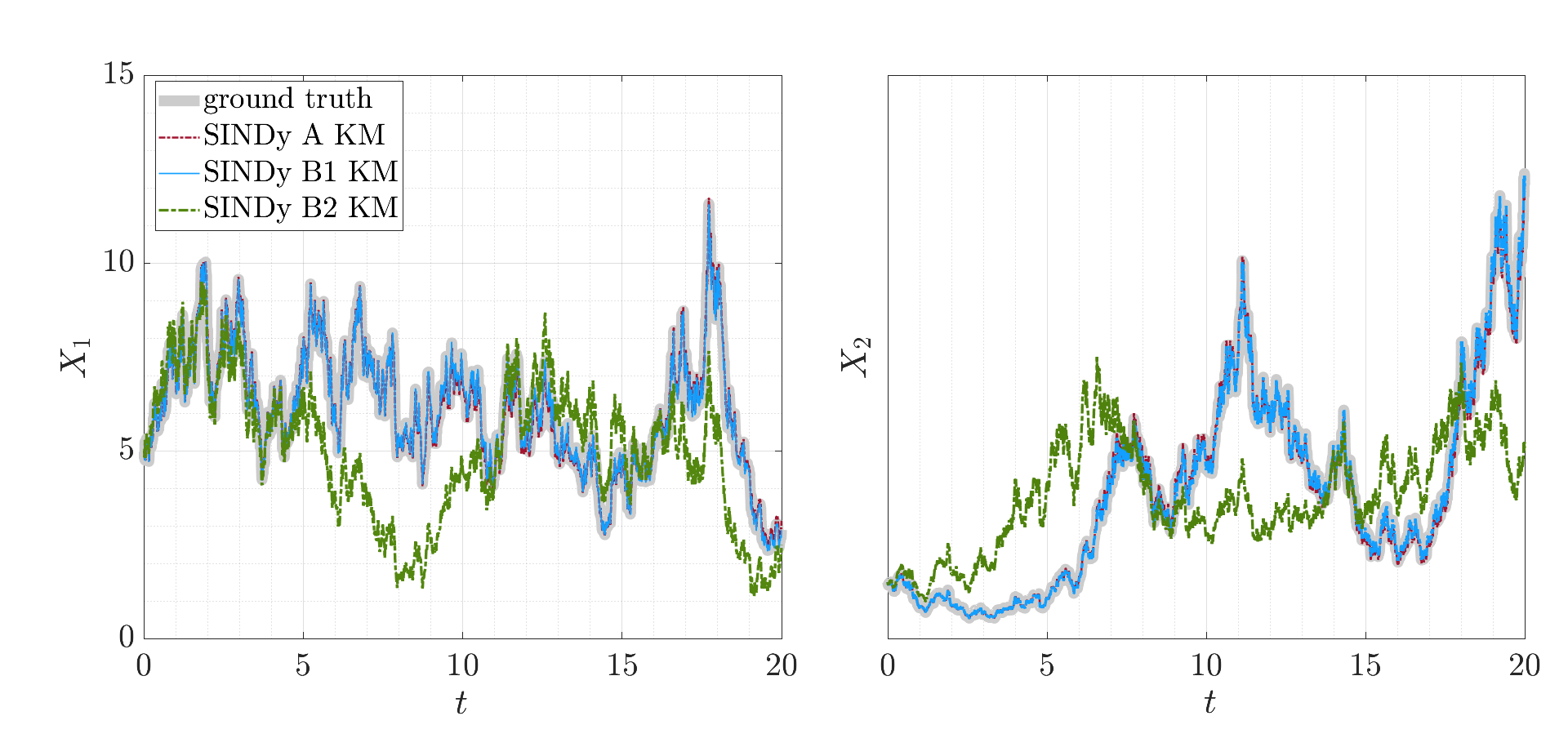}
    \caption{comparison of trajectories for the predator-prey SDDE \eqref{eq:pred_prey}: ground truth and approaches A, B1 and B2 all with KM.}\label{fig:pp_traj}
\end{figure}
    
\bigskip
Finally, as a third and last model we consider from \cite{kazmerchuk2007pricing} the SDDE
\begin{equation}\label{eq:pricing_option}
 \dd X(t)= rX(t)\dd t+\sigma(X(t),X(t-\tau))X(t)\dd W(t)
\end{equation}
modeling the price dynamics of a financial asset with memory, where $X(t)$ denotes the asset (or option) price at time $t$, $r=0.05$ is the constant risk-free interest rate, $\tau=0.002$, the function $\sigma$ represents the volatility (see below), which is a continuous function of also the history of the stock price and $W(t)$  is a standard Wiener process.
\eqref{eq:pricing_option} generalizes the classical Black–Scholes framework by incorporating historical dependence into the volatility term, capturing memory effects observed in real markets. Once more, data are obtained by simulating over the time window $[0,20]$ with time step $\Delta t=0.01$ via Euler–Maruyama, prescribing the constant history $X(s)\equiv100$. 80\% of the samples are used for training and 20\% for validation. For SINDy identification, custom libraries are constructed. The drift is identified using a polynomial library of degree $d=1$ and the diffusion is identified using a polynomial library of degree $d=4$. Both libraries are augmented with the logarithmic term $\ln(X(t) / X({t-\tau}))$ (see \eqref{eq:volatility} below) and $\lambda=0.025$ is used. Results of the full comparison are collected in Table \ref{tab:comparison3} and discussed in Section \ref{sec:discussion}.

As proposed in \cite{kazmerchuk2007pricing}, the volatility function $\sigma$ is based on a continuous-time analogue of the GARCH model, viz.
\begin{equation}\label{eq:volatility}
\sigma^{2}(X(t),X(t-\tau))= \frac{\gamma V}{\alpha + \gamma} + \frac{\alpha}{\tau (\alpha + \gamma)} \ln^{2}\left( \frac{X(t)}{X(t - \tau)} \right),    
\end{equation}
where $V=0.127$ is the long-run average variance, $\alpha=0.6$ and $\gamma=0.4$ are positive model parameters controlling the weighting of past return and reversion strength.
This formulation ensures that $\sigma$ remains strictly positive and reflects volatility clustering effects observed in empirical financial time series. Such memory-aware volatility plays a crucial role in capturing more realistic option price dynamics compared to models assuming constant or purely time-dependent volatility.
\begin{sidewaystable}[htbp]
\centering
\scriptsize
\begin{tabular}{ll c ccc ccc ccc ccc}
\toprule
& library & true & \multicolumn{12}{c}{$|\xi_{i}-{\bar{\xi_{i}}}|$} \\
\cmidrule(r){4-15} 
& terms & coeff. & \multicolumn{3}{c}{KM} & \multicolumn{3}{c}{FD} & \multicolumn{3}{c}{CD} & \multicolumn{3}{c}{TR} \\
\cmidrule(r){4-6} \cmidrule(r){7-9} \cmidrule(r){10-12} \cmidrule(r){13-15} 
& &  & $A$ & B1 & B2 & $A$ & B1 & B2 & $A$ & B1 & B2 & $A$ & B1 & B2  \\
\midrule
\multirow{1}{*}{drift} 
& $X(t)$ & $\xi_{2}=0.05$ & {$1.49e^{0}$} & \boldmath{$4.26e^{-6}$} & {$1.32e^{0}$} & {$1.48e^{0}$} & \boldmath{$3.65e^{-6}$} & {$1.32e^{0}$} & $3.03e^{-1}$ & \boldmath{$6.26e^{-6}$} & {$1.32e^{0}$}  & $3.09e^{0}$ & \boldmath{$7.57e^{-5}$} &  {$1.32e^{0}$} \\
\midrule
\multicolumn{15}{c}{\textbf{same diffusion}} \\
\midrule
diffusion
& $X^{2}(t)$ & $\xi_{5}=0.058$ & {$5.70e^{-1}$} & \boldmath{$2.35e^{-5}$} & {$9.20e^{-1}$} & {$5.70e^{-1}$} & \boldmath{$2.35e^{-5}$} & {$9.20e^{-1}$}  & {$5.70e^{-1}$} & \boldmath{$2.35e^{-5}$} &{$9.20e^{-1}$}  & {$5.70e^{-1}$} & \boldmath{$2.35e^{-5}$} & {$9.20e^{-1}$}  \\
& $G(Z(t))$ & $\xi_{33}=300$ & {$8.51e^{2}$} & \boldmath{$1.48e^{0}$} &  {$1.46e^{2}$} & {$8.51e^{1}$} & \boldmath{$1.48e^{0}$} & {$1.46e^{2}$} & {$8.51e^{1}$} & \boldmath{$1.48e^{0}$} & {$1.46e^{2}$}  & {$8.51e^{1}$} & \boldmath{$1.48e^{0}$} & {$1.46e^{2}$} \\
\midrule
\multicolumn{15}{c}{\textbf{performance metrics}} \\
\midrule
CPU [s] & \multicolumn{2}{c}{} & $\approx 4$ & $\approx 10 $ & $\approx 10$ & $\approx 4$ & $\approx 10$ & $\approx 10$ & $\approx 4$ & $\approx 10$ & $\approx 10$ & $\approx 4$ & $\approx 10 $ & $\approx 10$ \\
\midrule
& \multicolumn{2}{c}{full} & {$ 3.00e^{0} $} & \boldmath{$6.48e^{-6}$} & {$ 3.72e^{0} $}  & {$1.23e^{0}$} &  \boldmath{$9.67e^{-6}$} &  $4.37e^{0}$ & {$1.12e^{0}$} & \boldmath{$6.19e^{-6}$}  & $4.42e^{0}$ & {$7.88e^{0}$} & \boldmath{$6.46e^{-3}$}  & $2.39e^{1}$ \\
RMSE & \multicolumn{2}{c}{drift} & {$5.35e^{-1} $} &  \boldmath{$ 2.85e^{-5}$} & {$ 2.50e^{0} $} & {$7.72e^{-1}$} & \boldmath{$1.05e^{-5}$} & $2.32e^{0}$  & {$2.97e^{-1}$} & \boldmath{$7.83e^{-5}$} &  $2.48e^{0}$ & {$1.07e^{0}$} & \boldmath{$2.14e^{-4}$}  &  $9.68e^{1}$ \\
& \multicolumn{2}{c}{diffusion} & $9.00e^{0} $ & \boldmath{$1.17e^{-3}$} & {$ 2.67e^{0} $} & {$1.13e^{0}$} & \boldmath{$4.66e^{-3}$} & {$6.69e^{-1}$} & {$1.13e^{0}$} & \boldmath{$1.53e^{-3}$} & $1.29e^{0}$  & $9.00e^{0}$ & \boldmath{$2.11e^{-2}$}  & {$1.54e^{0}$} \\
\midrule
\multicolumn{15}{c}{\textbf{different diffusion}} \\
\midrule
diffusion
& $X^{2}(t)$ & $\xi_{5}=0.058$ & {$5.70e^{-1}$} & \boldmath{$2.35e^{-5}$} & {$9.20e^{-1}$} & {$3.08e^{0}$} & \boldmath{$7.37e^{-7}$} & {$2.50e^{-1}$} & {$7.97e^{-1}$} & \boldmath{$3.75e^{-6}$} &  {$5.67e^{-1}$} & {$6.16e^{-1}$} & \boldmath{$1.54e^{-5}$} & {$4.82e^{0}$} \\
& $G(Z(t))$ & $\xi_{33}=300$ & {$8.51e^{2}$} & \boldmath{$1.48e^{0}$} &  {$1.46e^{2}$}  & {$2.93e^{2}$} & \boldmath{$7.59e^{-1}$} & {$7.46e^{1}$} & {$2.94e^{2}$} & \boldmath{$5.43e^{-1}$} & {$2.46e^{2}$} & {$8.85e^{2}$} & \boldmath{$1.40e^{0}$} & {$2.90e^{2}$} \\
\midrule
\multicolumn{15}{c}{\textbf{performance metrics}} \\
\midrule
CPU [s] & \multicolumn{2}{c}{} & $\approx 4$ & $\approx 10 $ & $\approx 10$ & $\approx 4$ & $\approx 10$ & $\approx 10$ & $\approx 4$ & $\approx 10$ & $\approx 10$ & $\approx 4$ & $\approx 10 $ & $\approx 10$ \\
\midrule
& \multicolumn{2}{c}{full} & {$ 3.00e^{0} $} & \boldmath{$6.48e^{-6}$} & {$ 3.72e^{0} $} & {$3.67e^{1}$} &  \boldmath{$5.15e^{-6}$} & {$ 4.24e^{1} $}  & {$2.63e^{0}$} & \boldmath{$1.03e^{-6}$}  & {$ 2.91e^{0} $} & $7.76e^{0}$ & \boldmath{$4.40e^{-2}$}  & $2.76e^{1}$ \\
RMSE & \multicolumn{2}{c}{drift} & $5.35e^{-1} $ &  \boldmath{$ 2.85e^{-5}$} & {$ 2.50e^{0} $}  & { $5.34e^{-1}$} & \boldmath{$5.36e^{-4}$} &  {$ 5.88e^{0} $} & {$2.97e^{-1}$} & \boldmath{$5.27e^{-5}$} & {$ 2.28e^{0} $}  & {$1.07e^{0}$} & \boldmath{$7.83e^{-6}$}  & $9.13e^{1}$ \\
& \multicolumn{2}{c}{diffusion} & $9.00e^{0} $ & \boldmath{$1.17e^{-3}$} & {$ 2.67e^{0} $} & {$1.61e^{0}$} & \boldmath{$1.58e^{-5}$} & {$ 1.25e^{0} $} & {$6.47e^{1}$} & \boldmath{$1.88e^{-4}$} & {$ 4.43e^{0} $} & $9.30e^{0}$ & \boldmath{$4.15e^{-3}$}  & {$1.13e^{0}$} \\
\bottomrule
\end{tabular}
\caption{comparison of errors on estimated coefficients, RMSE on trajectories and CPU time for the option-pricing SDDE \eqref{eq:pricing_option} with \eqref{eq:volatility} for both drift and diffusion, all estimation methods KM, FD, CD and TR from Section \eqref{sec:ito_taylor} and all data reconstruction strategies A, B1 and B2 from Section \ref{sec:estimates}. Diffusion estimation might remain the same (KM in particular, first part above) or change (all of KM, FD, CD and TR, second part below). Recall that $G(Z(t))=X^{2}(t)\ln^{2}{(X(t)/X(t-\tau))}$.}
\label{tab:comparison3}
\end{sidewaystable}

\begin{figure}
    \centering
    \includegraphics[width=0.75\linewidth]{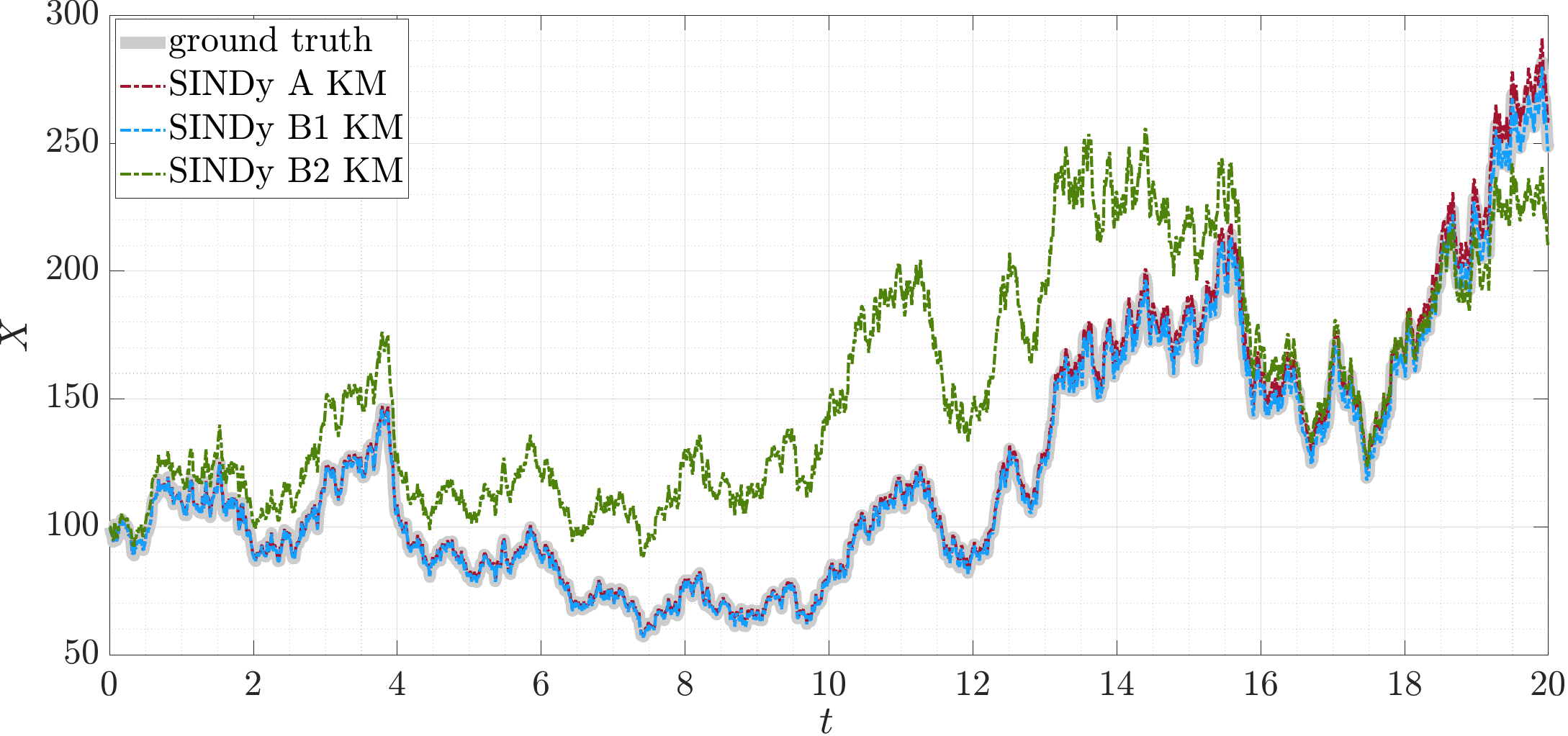}
    \caption{comparison of trajectories for the option pricing SDDE \eqref{eq:pricing_option}: ground truth and approaches A, B1 and B2 all with KM.}\label{fig:op_traj}
\end{figure}
\subsection{Discussion}\label{sec:discussion}
The analysis of the results collected in Tables \ref{tab:comparison1}, \ref{tab:comparison2} and \ref{tab:comparison3} reveals that the optimal estimation strategy is highly dependent on the specific model, the performance metric of interest (error on the reconstructed coefficients, RMSE on trajectories or CPU time) and the chosen methodology (A, B1 or B2 with KM, FD, CD or TR). In general the pre-regression averaging strategy of B1 overwhelmingly outperforms the other choices in terms of accuracy for all the three considered models. On the contrary, B2 consistently yields the highest absolute errors for every estimated coefficient and the highest RMSE values for the full system as well as for separated drift and diffusion and given the poor results, it is not considered a viable method. As far as A is concerned, the potential accuracy seems to depend on the prior available knowledge concerning the data via the (one-step) generation of relevant realizations. On the other hand, B1 operates without structural constraints, it frequently demonstrates a superior ability to identify the precise values of individual coefficients, thus highlighting its power in pure parameter and model discovery. Figures \ref{fig:log_traj} to \ref{fig:op_traj} present a comparison between the ground truth and the SINDy-predicted trajectories for each model. Figure \ref{fig:log_traj} (left) compares different estimation approaches for fixed KM, while Figure \ref{fig:log_traj} (right) compares various methods for estimation, holding method B1 constant. Based on the similar results for the other models, the KM approach was used for all subsequent figures (\ref{fig:pp_traj} and \ref{fig:op_traj}).

Figures \ref{fig:rmse_comp} through \ref{fig:log_comp} further clarify the comparative performance of the estimation strategies and numerical methods across different scenarios. Figure \ref{fig:rmse_comp} presents a comprehensive comparison of the RMSE for the state reconstructions across all methodologies as functions of the estimator (KM, FD, CD, TR) and data reconstruction approaches (A, B1, B2). It confirms the dominant accuracy of the pre-regression averaging approach B1 regardless of the estimator choice, with notably lower RMSE values reflecting superior reconstruction performance. Conversely, B2 again consistently yields higher RMSE values, reinforcing its unsuitability for reliable system identification. Figures \ref{fig:epsilon_variation} and \ref{fig:log_comp} in particular explore critical methodological factors influencing the accuracy and computational cost, further guiding practical implementation. Figure \ref{fig:epsilon_variation} examines the effect of the neighborhood tolerance parameter $\epsilon$ in B1 on the absolute errors of estimated coefficients and CPU time. This insight underlines the necessity of balancing statistical variance and bias during neighborhood selection to optimize performance. We anyway observe that decreasing $\epsilon$ too much does not lead to significant improvements and better results are obtained for diffusion rather than drift. Figure \ref{fig:log_comp} evaluates the impact of ensemble size M on estimation accuracy restricted only to the logistic SDDE model \eqref{eq:stoch_delay_logistic} (results are similar for the other models). We observe that, as the size M increases, Approach A does not notably improve the estimation of the diffusion coefficient; however, it does lead to better results for the drift term, but only for substantial increase of $M$, underlying the data-hungry nature of this approach. Conversely, Approach B2 consistently performs poorly across all values of M. In contrast, Approach B1 consistently outperforms the others in all cases, although its computational cost grows with increasing $M$ (indeed larger values of $M$ are not even considered).

Collectively, these figures validate the numerical superiority of B1 for SDDE model discovery. They provide actionable guidelines for tuning hyperparameters such as neighborhood size and ensemble cardinality, ensuring practical and effective application of the proposed SINDy framework for stochastic delay differential equations. These insights align with and reinforce the overall recommendation of prioritizing B1 in real-world discovery scenarios where model structures are unknown.

In conclusion, the selection of an optimal strategy depends entirely on the analytical objective. If the underlying model structure is known in advance and the goal is a rapid, computationally inexpensive estimation of the overall system behavior, A is a suitable choice. However, in a more realistic research scenario where the functional forms of the drift and diffusion terms are unknown and must be determined from the data, B1 is the decisively better option. Its ability to perform system identification without prior assumptions makes it a more robust, powerful and scientifically sound method for discovery-oriented analysis. The higher computational cost is justified by its capacity to derive the model directly from observations, minimizing user-introduced bias and leading to more accurate and reliable parameter estimation.
\begin{figure}
    \centering
    \includegraphics[width=1\linewidth]{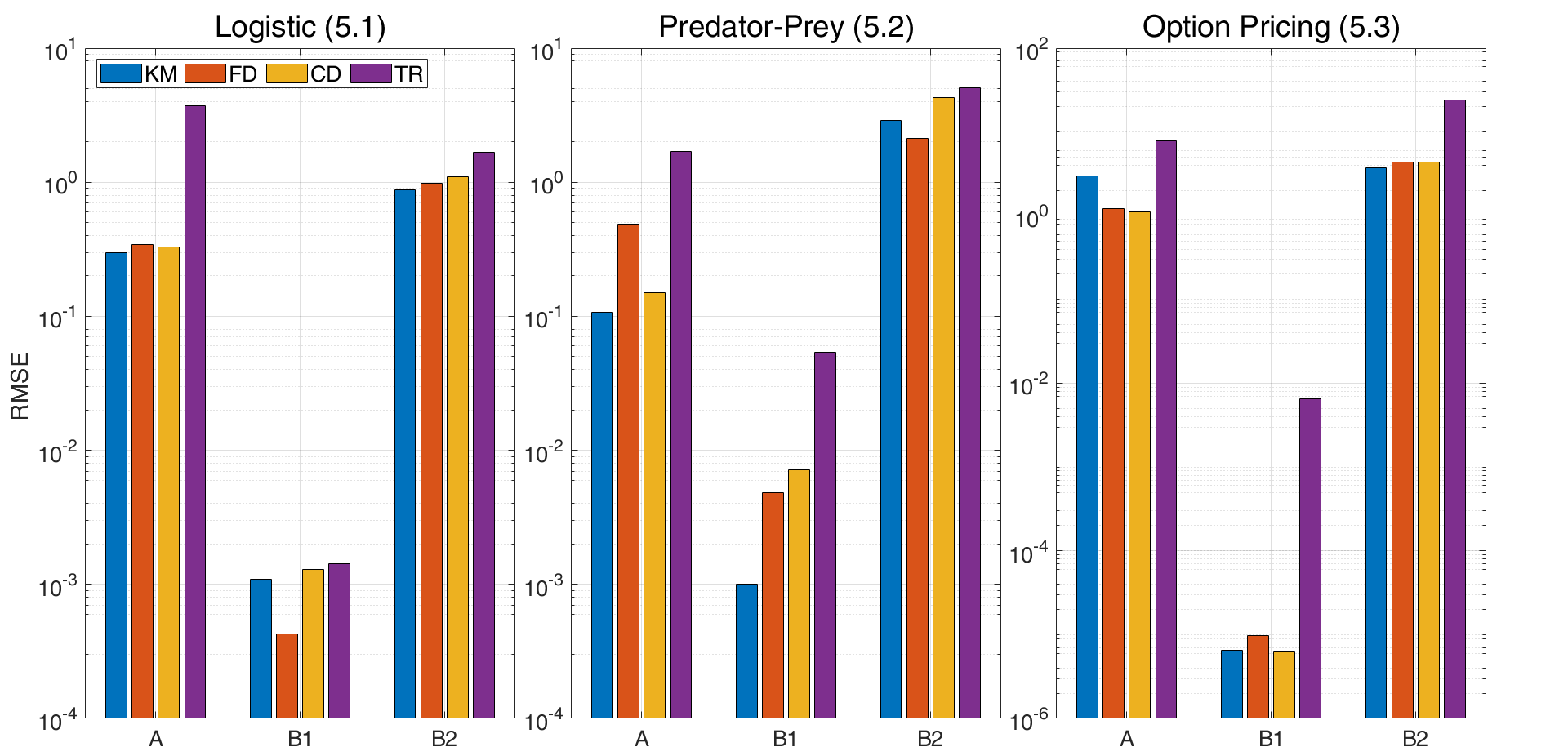}
    \caption{comparison of RMSE for the SDDEs \eqref{eq:stoch_delay_logistic}, \eqref{eq:pred_prey} and \eqref{eq:pricing_option} with all the approaches A, B1 and B2 against all the estimators KM, FD, CD and TR.}\label{fig:rmse_comp}
\end{figure}
\begin{figure}
    \centering
    \includegraphics[width=1\linewidth]{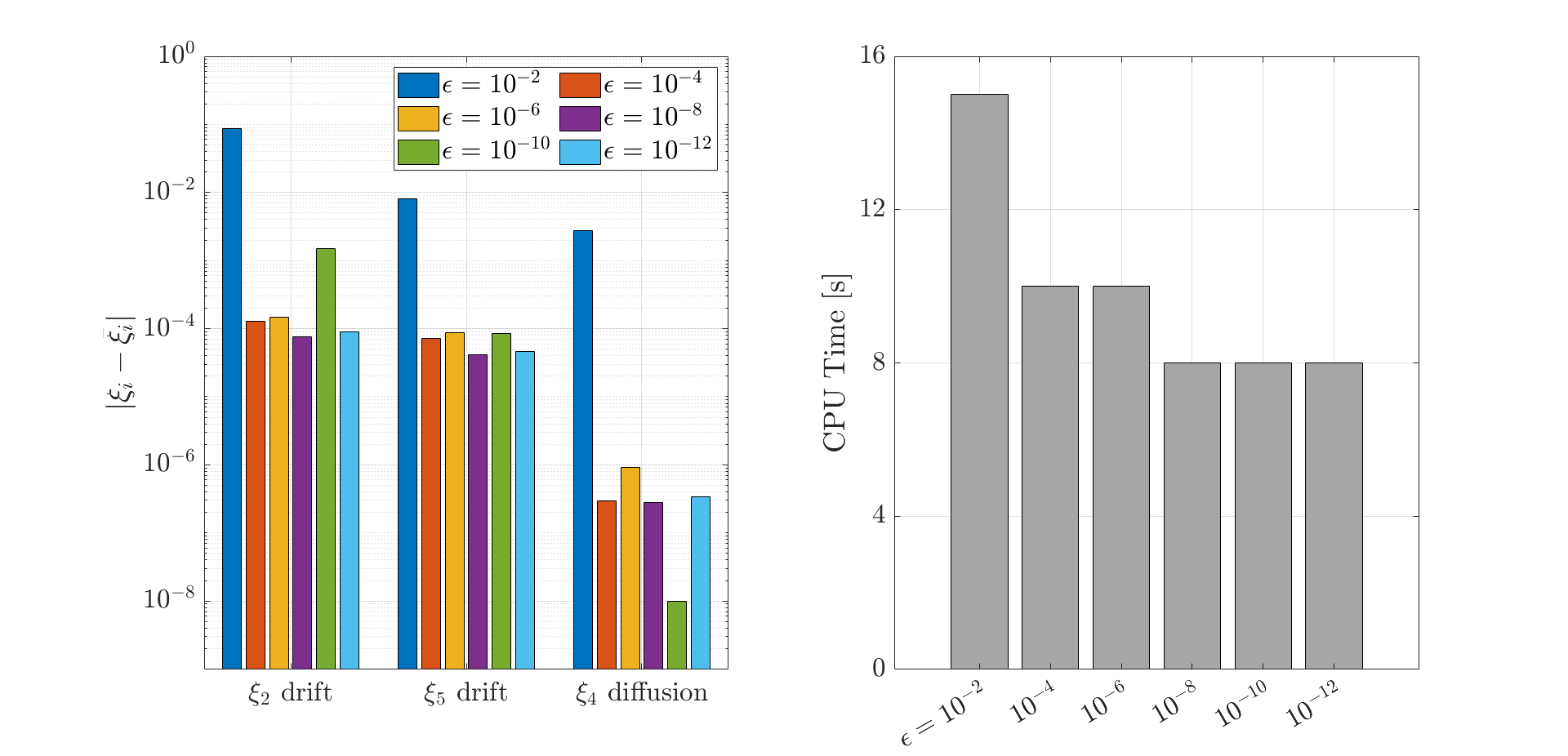}
    \caption{comparison of the absolute errors $|\xi_{i}-{\bar{\xi_{i}}}|$ (left) and CPU time (right) for varying $\epsilon$ for the logistic SDDE \eqref{eq:stoch_delay_logistic} with B1-FD.}\label{fig:epsilon_variation}
\end{figure}
\begin{figure}
    \centering
    \includegraphics[width=1\linewidth]{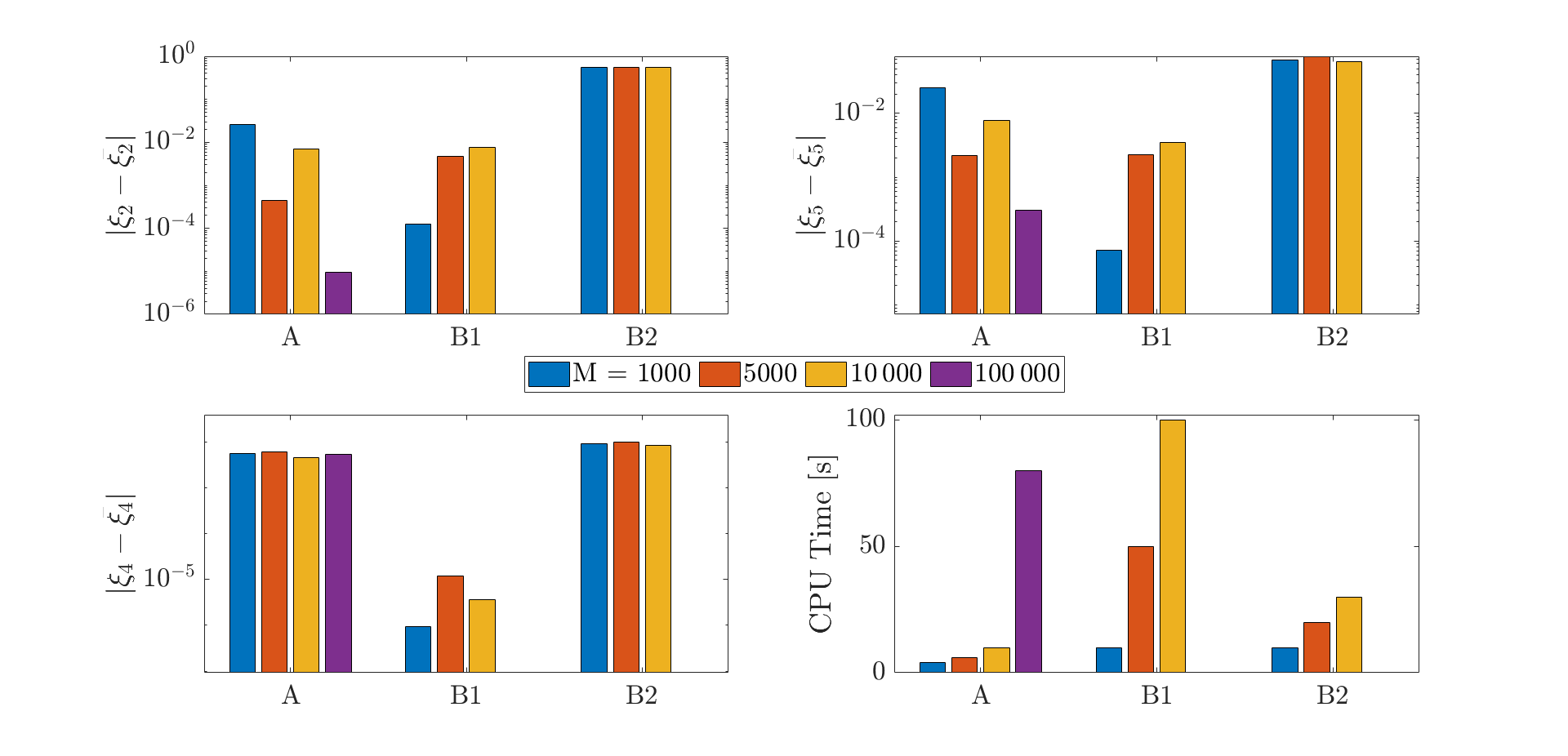}
    \caption{comparison of the absolute errors $|\xi_{i}-{\bar{\xi_{i}}}|$ for drift (top row) and diffusion (bottom left) for increasing $M$ for the logistic SDDE \eqref{eq:stoch_delay_logistic} using A, B1 and B2 all with FD.}\label{fig:log_comp}
\end{figure}
\section{Conclusion and future Work}\label{sec:conclusion}
This study presents a comprehensive evaluation of sparse identification strategies for SDDEs, focusing on the combined influence of estimator accuracy, trajectory ensemble construction and underlying system complexity. The comparative results underscore the superiority of exploiting the availability of multiple independent trajectories with pre-regression averaging (B1). Across all considered models including biological and financial systems with nonlinear delays this combination consistently yielded the lowest reconstruction errors and the most robust recovery of drift and diffusion dynamics. In contrast, strategies relying on single-path realizations-based reconstruction (A) or post-regression averaging (B2) performed reliably only under constrained or idealized conditions and often failed in more realistic, noise-sensitive scenarios. The analysis further reveals that estimator choice is not merely a computational convenience but a critical element influencing bias, robustness and interpretability. Together, these findings offer practical guidance for researchers applying sparse learning methods to complex stochastic systems with memory, highlighting the importance of ensemble-based strategies and carefully selected numerical estimators for successful model discovery and scientific insight.

Looking forward, the methodology offers promising opportunities for extension, including the incorporation of more complex and problem-specific candidate libraries to capture richer nonlinear and nonlocal dynamics, application and validation on real-world empirical datasets to address practical challenges such as measurement noise and missing data and the development of hybrid estimation techniques that combine the strengths of statistical and machine learning approaches. Additionally, integrating structure-preserving machine learning frameworks, such as physics-informed neural networks, could further enhance model interpretability and stability by respecting underlying physical constraints. Altogether, these avenues hold the potential to broaden the applicability and deepen the impact of sparse identification methods for SDDEs, facilitating robust, interpretable and computationally feasible data-driven discovery in complex stochastic systems with memory.

\section*{Acknowledgements}
DB, DC, RDA and MT are members of INdAM research group GNCS; DB is a
member of UMI research group ``Modellistica socio-epidemiologica''; DC and RDA are members of UMI research group ``Matematica per l'intelligenza Aartificiale e il machine learning''. The work of DB, DC, RDA and IS was supported by the European Union - NextGenerationEU with the project ``MOdellistica Numerica e Data-driven per l’Innovazione sostenibile - MONDI'' (CUP: G25F21003390007). The work of MT was supported by the Italian Ministry of University and Research (MUR) through a PhD grant PNRR DM351/22.}

\end{document}

\endinput